\documentclass{article}

\usepackage{PRIMEarxiv}

\usepackage[utf8]{inputenc} 
\usepackage[T1]{fontenc}    
\usepackage{hyperref}       
\usepackage{url}            
\usepackage{booktabs}       
\usepackage{amsfonts}       
\usepackage{nicefrac}       
\usepackage{microtype}      
\usepackage{lipsum}
\usepackage{fancyhdr}       
\usepackage{graphicx}       
\graphicspath{{media/}}     
\usepackage{siunitx}
\usepackage{amsfonts}
\usepackage{amsmath}
\usepackage{amsthm}
\usepackage{enumitem} 
\usepackage{url}
\usepackage{amssymb}

\newtheorem{theorem}{Theorem}

\newtheorem{corollary}{Corollary}

\newtheorem{definition}{Definition}

\newtheorem{lemma}{Lemma}

\newtheorem{remark}{Remark}

\title{Some results on Generalized HukuharaDiamond-$\alpha$ Derivative and Integral of Fuzzy Valued Functions and on Time Scales
}

\author{
 Selami BAYEĞ \\
  Department of Industrial Engineering\\
  University of Turkish Aeronautical Association \\
  Ankara, TÜRKİYE\\
  \texttt{sbayeg@thk.edu.tr} \\
   \AND
  Funda Raziye MERT \\
  Department of Software Engineering \\
  Adana Alparslan Türkeş Science and Technology University \\
  Adana, TÜRKİYE\\
  \texttt{rmert@atu.edu.tr} \\
   \And
   Billur KAYMAKÇALAN \\
   Department of Computer Engineering\\
   University of Turkish Aeronautical Association \\
  Ankara, TÜRKİYE\\
  \texttt{bkaymakcalan@thk.edu.tr} \\
}

\begin{document}
\maketitle

\begin{abstract}
In this paper, we define generalized Hukuhara diamond-alpha integral for fuzzy functions on
time scales and obtain some of its fundamental properties and also we establish
the relationship between diamond-alpha differentiation and integration.
\end{abstract}


\section{Introduction}

The theory of dynamic equations on time scales is a relatively recent field of research that has witnessed significant growth over the past 35 years. Time scale theory was developed to bridge continuous and discrete structures, enabling the simultaneous treatment of both difference and differential equations. This theory expands upon traditional differential equations to include dynamic equations, thereby providing a unified framework. An overview of the basics of time scale calculus and some notable recent developments can be found in \cite{Hilger1, Hilger2, Bohner, Bohner2, Agarwal, Guseinov, Guseinov2, Anderson, Malinowska, Malinowska2, Lakshmikantham2, Kayar}.

In many real-world phenomena, various uncertainties must be considered to achieve an accurate understanding. To address this, Zadeh \cite{Zadeh} introduced fuzzy set theory, which allows for the representation of vague or imprecise concepts. The theory of fuzzy differential equations (FDEs) was subsequently established by Kaleva \cite{Kaleva}, and by Lakshmikantham and Mohapatra \cite{Lakshmikantham}, with numerous applications explored thereafter.

A limitation of methods based on Hukuhara differentiability is that the solutions to FDEs are typically valid only for longer support intervals. To overcome this, Bede et al. \cite{Bede2} studied generalized Hukuhara differentiability, which has attracted the attention of many researchers \cite{Bede3, Li, Stefanini2} due to its advantages in handling fuzzy number-valued functions. 

In the context of fuzzy dynamic equations on time scales, Vasavi et al. \cite{Vasavi, Vasavi1, Vasavi2, Vasavi3} introduced the Hukuhara, second-type Hukuhara, and generalized delta derivatives by implementing the Hukuhara difference. However, these derivatives are limited to fuzzy number-valued functions on time scales where the diameter increases with length. Fard and Bidgoli \cite{Fard} further investigated the calculus of fuzzy functions on time scales.

Delta and nabla derivatives have typically been used independently to analyze fuzzy number-valued functions on time scales. You et al. introduced the nabla Hukuhara derivative of fuzzy functions on time scales and showed some of its basic properties in \cite{you}. Mert et al. introduced generalized Hukuhara delta derivative of fuzzy functions on time scales and they provide several new characterizations of generalized Hukuhara delta differentiable fuzzy functions on time scales, derived from the delta differentiability of their endpoint functions in \cite{mert1}. 

In certain cases---especially those involving time scales with discrete points---these derivatives may only describe changes in functions on one side of a point. To address this limitation, dynamic derivatives known as the \textit{diamond-alpha derivative} were introduced by Sheng et al. \cite{Sheng}. This derivative is a convex linear combination of the delta and nabla derivatives. Later, Rogers et al. \cite{Rogers} redefined the diamond-alpha derivative independently of the standard delta and nabla derivatives and further investigated its properties. In \cite{Truong}, the diamond-alpha derivative was extended to interval-valued functions on time scales via generalized Hukuhara differences. That study also investigated a particular class of interval differential equations related to the diamond-alpha derivative.

In \cite{bayeg}, we introduced the $\diamond_{gH}^{\alpha}$-derivative for fuzzy number-valued functions defined on time scales and explored its properties under various conditions imposed on the time scale $\mathbb{T}$. Based on this work, some of the main results presented in this paper include the following: We demonstrated that the existence of the $\Delta_{gH}$-derivative of a fuzzy function at a given point does not necessarily imply the existence of the $\nabla_{gH}$-derivative at that point. Moreover, we established that if a fuzzy function possesses a $\diamond_{gH}^{\alpha}$-derivative at a point, then it must be continuous at that point. Additionally, we proved that $\diamond_{gH}^{\alpha}$-differentiability implies both $\Delta_{gH}$ and $\nabla_{gH}$ differentiability. Furthermore, we derived essential results such as the addition, scalar multiplication, and product rules for $\diamond_{gH}^{\alpha}$-differentiability.

The study of integrals of fuzzy functions on time scales has gained significant attention in recent years, primarily due to the need for more generalized mathematical frameworks capable of handling uncertainty and variability in dynamic systems. One foundational contribution in this area is the introduction of the generalized delta integral for fuzzy set-valued functions, which enables the integration of fuzzy functions defined on time scales while addressing properties such as normality, convexity, and upper semicontinuity of the fuzzy sets involved \cite{Vasavi1, Vasavi2, Vasavi3}. Furthermore, the fuzzy C-delta integral was developed to extend the concept of integrability to interval-valued and fuzzy-valued functions on time scales, establishing necessary and sufficient conditions for such integrability \cite{You2017}. Using the generalized Hukuhara difference, Fard and Bidgoli recently introduced and investigated Henstock--Kurzweil integrals of fuzzy-valued functions on time scales \cite{Fard}.

The nabla integral, another important construct, has been defined for fuzzy functions, highlighting the relationship between nabla differentiation and integration and thereby enriching the calculus of fuzzy functions on time scales \cite{Leelavathi2018, Leelavathi}. Moreover, the fuzzy Henstock--Kurzweil--Stieltjes integral has been introduced, incorporating properties of fuzzy functions within the context of time scales and allowing for the application of convergence theorems, such as the monotone and dominated convergence theorems \cite{Afariogun2021, Li2024}. This integral framework is particularly useful for addressing complex systems where traditional integration methods may fall short.

The exploration of fuzzy Volterra integral equations on time scales further illustrates the versatility of these integrals in modeling dynamic systems. These studies have provided existence theorems and established equivalences to linear fuzzy dynamic equations \cite{Shahidi2020, Shahidi2023}. Overall, the integration of fuzzy functions on time scales not only broadens the scope of mathematical analysis but also enhances the applicability of fuzzy calculus across various fields, including control theory and dynamic modeling.

In this work, we introduce the generalized Hukuhara diamond-alpha integrals for fuzzy functions on time scales and examine their fundamental properties, such as the Fundamental Theorem of Calculus and the method of Integration by Parts, motivated by \cite{Vasavi1, Vasavi2, Vasavi3, Fard, you, mert1, Leelavathi2018, Leelavathi,mert2}.

This paper is structured as follows: Section~2 provides foundational concepts related to fuzzy sets and time scales. Section~3 presents key results concerning the generalized Hukuhara diamond-alpha derivative. In Section~4, we introduce the definition of the generalized Hukuhara diamond-alpha integral, explore its fundamental properties, and discuss its connection to the Fundamental Theorem of Calculus and the method of Integration by Parts. Finally, Section~5 concludes the paper with a summary of findings and potential directions for future research.

\section{Preliminaries}

\begin{definition}\label{def1}\cite{Bohner}
	A non-empty closed subset of the real numbers $\mathbb R$ is called a time scale, often denoted by $\mathbb T$.
	
\end{definition}

\begin{definition}\label{def2a}\cite{Bohner}
	The function \( \sigma : \mathbb{T} \rightarrow \mathbb{R} \) defined by
	\begin{equation*}
		\sigma(t) = \inf\{s \in \mathbb{T} : s > t \}
	\end{equation*}
	is called the forward jump operator. Additionally, we set \( \inf\varnothing := \sup\mathbb{T} \).
\end{definition}

\begin{definition}\label{def2b}\cite{Bohner}
	The function \( \rho : \mathbb{T} \rightarrow \mathbb{R} \) defined by
	\begin{equation*}
		\rho(t) = \sup\{s \in \mathbb{T} : s < t \}
	\end{equation*}
	is called the backward jump operator. Additionally, we set \( \sup\varnothing := \inf \mathbb{T} \).
\end{definition}

 \begin{definition}\label{def2*}\cite{Bohner}
If \( \sigma(t) > t \), then \( t \in \mathbb{T} \) is said to be a  right-scattered point.
\end{definition}

 \begin{definition}\label{def2**}\cite{Bohner}
If \( \rho(t) < t \), then \( t \in \mathbb{T} \) is said to be a left-scattered point.
\end{definition}
 \begin{definition}\label{def2***}\cite{Bohner}
 If \( \sigma(t) = t \) and \( t \neq \sup\mathbb{T} \), then \( t \in \mathbb{T} \) is said to be a right-dense point.
\end{definition}
 \begin{definition}\label{def2***}\cite{Bohner}
If \( \rho(t) = t \) and \( t \neq \inf\mathbb{T} \), then \( t \in \mathbb{T} \) is said to be a left-dense point.
\end{definition}
\begin{definition}\cite{Bohner}
	The function $\mu :\mathbb{T}\rightarrow [0,\infty)$ defined by $\mu(t) = \sigma(t)-t$ is called the (forward) graininess.
\end{definition}
\begin{definition}\label{def3}\cite{Bohner}
	The function $\nu : \mathbb T\rightarrow [0,\infty)$ defined by $\nu(t)= t-\rho(t)$ is called the backward graininess.
\end{definition}
Additionally, we define the following notations for simplicity in the definitions and theorems throughout this paper:  $\mu _{st}=\sigma (s)-t$  and $\nu _{st}=t-\rho (s)$.\\ 
The set \( \mathbb{T}^{\kappa} \) is defined as follows: if \( \mathbb{T} \) has a left-scattered maximum \( m \), then \( \mathbb{T}^{\kappa} := \mathbb{T} \setminus \{m\} \). If no such maximum exists, then \( \mathbb{T}^{\kappa} := \mathbb{T} \). Similarly, the set \( \mathbb{T}_{\kappa} \) is defined as follows: if \( \mathbb{T} \) has a right-scattered minimum \( m \), then \( \mathbb{T}_{\kappa} := \mathbb{T} \setminus \{m\} \). If no such minimum exists, then \( \mathbb{T}_{\kappa} := \mathbb{T} \).

\begin{definition}\label{def4}\cite{Bohner}
	Let \( h : \mathbb{T} \rightarrow \mathbb{R} \) be a function and let \( s \in \mathbb{T}^{\kappa} \). We define \( (\Delta h)(s) \) as the number (if it exists) that satisfies the following property: for any \( \epsilon > 0 \), there exists a neighborhood \( N_{\mathbb{T}} \) of \( s \) given by \( N_{\mathbb{T}} := (s - \delta, s + \delta) \cap \mathbb{T} \) for some \( \delta > 0 \) such that
	$$|[h(\sigma(s)) - h(t)] - (\Delta h)(s)[\sigma(s) - t]| \leq \epsilon |\sigma(s) - t|$$
	for all \( t \in N_{\mathbb{T}} \). The value \( (\Delta h)(s)\) is called the delta derivative of \( h \) at \( s \).
\end{definition}

\begin{definition}\label{def5}\cite{Bohner}
	Let \( h : \mathbb{T} \rightarrow \mathbb{R} \) be a function, and let \( s \in \mathbb{T}_{\kappa} \). We define \( (\nabla h)(s) \) as the number (if it exists) that satisfies the following property: for any \( \epsilon > 0 \), there is a neighborhood given by \( N_{\mathbb{T}} := (s - \delta, s + \delta) \cap \mathbb{T} \) for some \( \delta > 0 \) such that
	$$|[h(\rho(s)) - h(t)] -(\nabla h)(s)[\rho(s) - t]| \leq \epsilon |\rho(s) - t|$$
	for all \( t \in N_{\mathbb{T}} \). The value \( (\nabla h)(s) \) is referred to as the nabla derivative of \( h \) at \( s \).
\end{definition}

\begin{definition} \label{def6}\cite{Rogers}
	Let \( h : \mathbb{T} \rightarrow \mathbb{R} \) be a function and \( s \in \mathbb{T}^{\kappa} \cap \mathbb{T}_{\kappa} \). Then the \( \diamond^{\alpha} \)-derivative of \( h \) at the point \( s \in \mathbb{T}_{\kappa}^{\kappa} \), denoted by \( (\diamond^{\alpha}h)(s) \), is the number (provided it exists) that satisfies the following property: for any \( \epsilon > 0 \), there is a neighborhood given by \( N_{\mathbb{T}} := (s - \delta, s + \delta) \cap \mathbb{T} \) for some \( \delta > 0 \) such that 
	$$\begin{aligned}
		| \alpha |h(\sigma(s)) - h(t)| |\nu_{st}| + (1 - \alpha)|h(\rho(s)) - h(t)| |\mu_{st}| -  (\diamond^{\alpha}h)(s) |\nu_{st} \mu_{st}| \leq \epsilon |\nu_{st} \mu_{st}|,
	\end{aligned}$$
for any \( t \in N_{\mathbb{T}} \). Here, \( (\diamond^{\alpha}h)(s)\) referred to as the diamond-alpha derivative of \( h \) at \( s \).
\end{definition}

\begin{definition}\label{def7} \cite{Bourdin} A time scale $\mathbb T$ is called quasi-regular if $\sigma (\rho (t))=t$
for all $t\in \mathbb T_{\kappa }$ and $\sigma (\rho (t))=t$ for all $t\in \mathbb T^{\kappa
}$.
\end{definition}

\begin{definition}\label{def8}
A time scale $\mathbb T$ is called a homogeneous if $\mu (t)=\nu (t)=c\geq 0$ for
all $t\in \mathbb T.$
\end{definition}

\begin{remark}\label{rem1}
Every homogeneous time scale is regular.
\end{remark}

\begin{theorem}\label{thm1}
Suppose $\mathbb T$ is a time scale and $f$ is a continuous function on $\mathbb T$. Then,

\begin{enumerate}
\item If $f$ is $\Delta $-differentiable on $\mathbb T^{\kappa }$, then $f$ is $%
\nabla $-differentiable and $(\nabla f)(t)=(\Delta f)(\rho (t))$ for all $%
t\in \mathbb T_{\kappa }$ such that $\sigma (\rho (t))=t$.

\item If $f$ is $\nabla $-differentiable on $\mathbb T_{\kappa }$, then $f$ is $%
\Delta $-differentiable and $(\Delta f)(t)=(\nabla f)(\sigma (t))$ for all $%
t\in \mathbb T^{\kappa }$ such that $\rho (\sigma (t))=t$.
\end{enumerate}
\end{theorem}

\begin{lemma}\label{lem1}
Suppose that $\mathbb T$ is a quasi-regular time scale and the functions $f$, $%
f^{\Delta }$, and $f^{\nabla }$ are continuous on $\mathbb T_{\kappa }^{\kappa }.$
Then, the followings hold:
\end{lemma}

\begin{enumerate}
\item $(\nabla \sigma )(t)$ $(\Delta \rho )(t)$ $=1$.

\item $\nu (t)\nabla ((\Delta f)(t))=\mu (t)\Delta ((\nabla f)(t)).$
\end{enumerate}

\begin{remark}\label{rem2}
Assume $\sigma $ is $\nabla $-differentiable on $\mathbb T_{\kappa }$ with $(\nabla
\sigma )(t)=m$. Since $\sigma $ is an increasing function, then $m\geq 1$
for $t\geq 0$ and $0<m\leq 1$ for $t\leq 0.$
\end{remark}

\begin{corollary}\label{corol1}
If a time scale $\mathbb T$ is homogeneous and the functions $f$, $\Delta f$, and $%
\nabla f$ are continuous on $\mathbb T_{\kappa }^{\kappa }$, then

\begin{enumerate}
\item $(\nabla \sigma )(t)$ $=1$ and $(\Delta \rho )(t)$ $=1$.

\item $\nabla ((\Delta f)(t))=\Delta ((\nabla f)(t)).$
\end{enumerate}
\end{corollary}

\begin{remark}\label{rem3*}
Assume $\sigma $ is $\nabla $-differentiable on $\mathbb T_{\kappa }$ with $(\nabla
\sigma )(t)$ $=m$. Since $\sigma $ is an increasing function, then $m\geq 1$
for $t\geq 0$ and $0<m\leq 1$ for $t\leq 0$.
\end{remark}

\begin{definition}\cite{Zadeh} \label{def9}
	A fuzzy set \( u \) in a universe of discourse \( U \) is represented by a function $u:U\rightarrow [0,1]$, where \( u(x) \) is the membership degree of \( x \) to the fuzzy set \( u \).
\end{definition}
We use \( \mathcal F(U) \) to denote the set of all fuzzy subsets of \( U \).
\begin{definition}\cite{Negoita}\label{def10}
	Let $u:U\rightarrow \lbrack 0,1]$ be a fuzzy set. The r-level sets of $u$
	are defined as
	\[
	u_{r}=\left\{ x\in U:u(x)\geq r\right\}
	\]%
	for $0<r\leq 1$. The $0$-level set of $u$ 
	\[
	u_{0}=cl\left\{ x\in U:u(x)>0\right\} 
	\]%
	is called the support of the fuzzy set $u$. Here, $cl$ denotes the closure
	of the set $u.$
\end{definition}

\begin{definition}\cite{Negoita}\label{def11}
	Let \( u: \mathbb{R} \rightarrow [0,1] \) be a fuzzy subset of the real numbers. Then, \( u \) is said to be a fuzzy number if it fulfills the following criteria:

	\begin{enumerate}
		\item \( u \) is normal, which means that there exists an \( x_{0} \in \mathbb{R} \) such that \( u(x_{0}) = 1 \).
		
		\item \( u \) is quasiconcave, which means that for all \( \lambda \in [0,1] \), the inequality \( u(\lambda x + (1-\lambda)y) \geq \min\{u(x), u(y)\} \) holds.
		
		\item \( u \) is upper semicontinuous on \( \mathbb{R} \), which means that for any \( \epsilon > 0 \), there exists a \( \delta > 0 \) such that \( u(x) - u(x_{0}) < \epsilon \) whenever \( |x - x_{0}| < \delta \).
		
		\item \( u \) is compactly supported, which means that the closure \( \text{cl}\{ x \in \mathbb{R} : u(x) > 0 \} \) is compact.
	\end{enumerate}
\end{definition}
We use \( \mathcal F_N(\mathbb{R}) \) to denote the set of all fuzzy numbers of \( \mathbb{R} \).

\begin{definition}\label{def12}
	Let $a_1\leq a_2\leq a_3$ be real numbers. The fuzzy number denoted by $u=(a_1,a_2,a_3)$ is called a triangular fuzzy number whose membership function is
	\[
	u(x)=\left\{ 
	\begin{array}{cc}
		\frac{x-a_1}{a_2-a_1}, & a_1\leq x\leq a_2, \\ 
		\frac{a_3-x}{a_3-a_2}, & a_2\leq x\leq a_3, \\ 
		0, & otherwise.%
	\end{array}%
	\right. 
	\]
\end{definition}

\begin{definition}\cite{Stefanini}\label{def13}
	Let $u,v\in $ $\mathcal F_N(\mathbb{R})$. The generalized Hukuhara difference (gH-difference) is the fuzzy number $w$, if it exists, such that%
	$$
	u\ominus _{gH}v=w\Longleftrightarrow u=v+w \text{   or   } v=u+(-1)w. $$

\end{definition}
Since level sets of a fuzzy number are closed and bounded intervals, we will denote r-level set of a fuzzy number $u$ by $u_{r}=[u_{r}^{-},u_{r}^{+}]$ and its length by $len(u_{r})=u_{r}^{+}-u_{r}^{-}$.

\begin{remark}\label{rem3}
    The criteria for the existence of \( w = u \ominus_{g_H} v \) in \( \mathcal F_N(\mathbb{R}) \) are as follows:

    \textbf{Case (i):}  
    \begin{itemize}
        \item \( w_r^- = u_r^- - v_r^- \) and \( w_r^+ = u_r^+ - v_r^+ \)
        \item Here, \( w_r^- \) must be increasing, \( w_r^+ \) must be decreasing, and it must hold that \( w_r^- \leq w_r^+ \) for all \( r \) in \([0, 1]\).
    \end{itemize}

    \textbf{Case (ii):}  
    \begin{itemize}
        \item \( w_r^- = u_r^+ - v_r^+ \) and \( w_r^+ = u_r^- - v_r^- \)
        \item Similarly, \( w_r^- \) should be increasing, \( w_r^+ \) should be decreasing, and \( w_r^- \leq w_r^+ \) must hold for all \( r \) in \([0, 1]\).
    \end{itemize}
\end{remark}

\begin{theorem}\cite{Stefanini, Bede}\label{thm3} Let $u,v\in $ $\mathcal F_N(%
	\mathbb{R}
	).$ If gH-difference $u\ominus _{gH}v\in $ $\mathcal F_N(%
	\mathbb{R}
	)$ exists, then 
	\[
	\left( u\ominus _{gH}v\right) _{r}=[\min
	\{u_{r}^{-}-v_{r}^{-},u_{r}^{+}-v_{r}^{+}\},\max
	\{u_{r}^{-}-v_{r}^{-},u_{r}^{+}-v_{r}^{+}\}].
	\]
\end{theorem}

\begin{definition}\cite{Diamond2} \label{def14}
	The metric $D_{\infty }:\mathcal F_N(%
	\mathbb{R}
	)\times \mathcal F_N(%
	\mathbb{R}
	)\rightarrow 
	\mathbb{R}
	^{+}\cup \{0\}$ defined by
	\[
	D_{\infty }(u,v)=\sup_{r\in \lbrack 0,1]}D(u_{r},v_{r}), 
	\]
	where $D(u_{r},v_{r})=\max \left\{ \left\vert
	u_{r}^{-}-v_{r}^{-}\right\vert ,\left\vert u_{r}^{+}-v_{r}^{+}\right\vert
	\right\}$, $u_{r}=[u_{r}^{-},u_{r}^{+}]$, $v_{r}=[v_{r}^{-},v_{r}^{+}]$, is called Hausdorff metric for fuzzy numbers.
\end{definition}
The Hausdorff metric provides a way to measure the distance between two fuzzy sets by considering their level sets. This metric allows researchers to compare the similarity or dissimilarity of fuzzy sets in a rigorous mathematical way. Specifically, it can be used to quantify how far apart two fuzzy sets are based on their support and their membership functions.
\begin{theorem}\cite{Diamond2} \label{prop3}

	Let $a$, $b$, $c$, $d$ $\in \mathcal F_N(\mathbb{R})$ and $m\in\mathbb{R}$. The Hausdorff metric satisfies the followings:
	\begin{enumerate}
		\item $D_{\infty }\left( a+c,b+c\right) =D_{\infty }\left( a,b\right).$

  	\item $D_{\infty }\left( ma,mb\right) =\left\vert m\right\vert
		D_{\infty }\left( a,b\right).$

		\item $D_{\infty }\left( a+b,c+d\right) \leq D_{\infty }\left(
		a,c\right) +D_{\infty }\left( b,d\right). $
	\end{enumerate}
	
\end{theorem}

\begin{definition} \cite{ mert1}\label{def13}
	Let $f:\mathbb{T}\rightarrow F(
	\mathbb{R}
)$ be a fuzzy function and $s\in \mathbb{T}^{\kappa }$. The generalized Hukuhara
	delta derivative of $f$ \ at $s$, if it exists, is a fuzzy number $\left(
	\Delta _{gH}f\right) (s)\in  (
	\mathbb{R}
)$ such that for any given $\epsilon >0,$ there exists a neighborhood $
	N_{\mathbb{T}}(s,\delta )=(s-\delta ,s+\delta )\cap \mathbb{T}$ for some $\delta >0,$ such
	that for all $t\in N_{\mathbb{T}}$,  $f(\sigma (s))\ominus _{gH}f(t)$ exists and  we
	have%
	\[
	D_{\infty }(f(\sigma (s))\ominus _{gH}f(t),(\Delta _{gH}f)(s)\mu
	_{st})\leq \epsilon \left\vert \mu _{st}\right\vert .
	\]

\end{definition}

\begin{definition}\cite{you}\label{def14}
	Let $f:\mathbb{T}\rightarrow F(
	\mathbb{R})$ be a fuzzy function and $s\in \mathbb{T}_{\kappa }$. The generalized Hukuhara
	nabla derivative of $f$ \ at $s$, if it exists, is a fuzzy number $\left(
	\nabla _{gH}f\right) (s)\in (\mathbb{R})$ such that for any given $\epsilon >0,$ there exists a neighborhood $%
	N_{\mathbb{T}}(s,\delta )=(s-\delta ,s+\delta )\cap \mathbb{T}$ for some $\delta >0,$ such
	that for all $t\in N_{\mathbb{T}}$,  $f(t)\ominus _{gH}f(\rho (s))$ exists and we have
	\[
	D_{\infty }(f(t)\ominus _{gH}f(\rho (s)),(\nabla _{gH}f)(s)\nu _{st})\leq
	\epsilon \left\vert \nu _{st}\right\vert .
	\]
\end{definition}

\begin{definition}\cite{bayeg}\label{def15}
	Let $f:\mathbb{T}\rightarrow F(
	\mathbb{R}	)$ be a fuzzy function and $s\in \mathbb{T}_{\kappa }^{\kappa }$. The generalized
	Hukuhara diamond-$\alpha $ derivative of $f$ \ at $s$, if it exists, is a
	fuzzy number $\left( \diamond _{gH}^{\alpha }f\right) (s)\in (%
	\mathbb{R}
	)$ such that for any given $\epsilon >0,$ there exists a neighborhood $%
	N_{\mathbb{T}}(s,\delta )=(s-\delta ,s+\delta )\cap \mathbb{T}$ for some $\delta >0,$ such
	that for all $t\in N_{\mathbb{T}}$  $f(\sigma (s))\ominus _{gH}f(t)$ and  $%
	f(t)\ominus _{gH}f(\rho (s))$ exist and we have 
	\[
	D_{\infty }(\alpha \lbrack f(\sigma (s))\ominus _{gH}f(t)]\nu
	_{st}+(1-\alpha )[f(t)\ominus _{gH}f(\rho (s))]\mu _{st},(\diamond
	_{gH}^{\alpha }f)(s)\mu _{st}\nu _{st})\leq \epsilon \left\vert \mu _{st}\nu
	_{st}\right\vert .
	\]
\end{definition}

\section{Some Results on Generalized Hukuhara Diamond Alpha Derivative}

\begin{theorem}\label{thm4}
Let  $f:\mathbb T\rightarrow \mathcal F_N(\mathbb R)$ be a fuzzy function on a time scale $\mathbb T$. The
existence of the $\nabla _{gH}$-derivative of $f$ at $s\in \mathbb T_{\kappa
}^{\kappa }$ does not imply the existence of the $\Delta_{gH}$-derivative of $f$
at $s\in \mathbb T_{\kappa }^{\kappa }$ , and vice versa.

\begin{proof}
Let us consider the function%
\begin{equation*}
f(t)=\left\{ 
\begin{array}{cc}
(1,2,3)t\sin (1/t), & t\neq 0, \\ 
0, & t=0%
\end{array}%
\right. 
\end{equation*}%
on $\mathbb T=[-3,-1]\cup \lbrack 0,2]$. The $r$-level sets of $f$ are%
\begin{equation*}
f_{r}(t)=\left\{ 
\begin{array}{cc}
\lbrack 1+r,3-r]t\sin (1/t), & t\neq 0, \\ 
0, & t=0.%
\end{array}%
\right. 
\end{equation*}%
The function $f$ is continuous at $\ 0$ and the point $0\in \mathbb T$ is
right-dense and left-scattered. By \cite{Leelavathi}, $f$ is $\nabla _{gH}$-differentiable
at $0$. But the limit 
\begin{eqnarray*}
\lim_{t\rightarrow s}{\frac{f(\sigma (s))\ominus_{gH}f(t)}{\sigma(s)-t}}
&=&\lim_{t\rightarrow 0}\frac{f(\sigma(0))\ominus_{gH}f(t)}{\sigma(0)-t} \\
&=&\lim_{t\rightarrow 0^{+}}\frac{f(0)\ominus _{gH}f(t)}{0-t}
\end{eqnarray*}%
does not exist. Since, $f(0)\ominus _{gH}f(t)=h(t)$ $\Leftrightarrow
h_{r}^{-}(t)=\min
\{f_{r}^{-}(0)-f_{r}^{-}(t),f_{r}^{+}(0)-f_{r}^{+}(t)\}=\min
\{-f_{r}^{-}(t),-f_{r}^{+}(t)\}=-f_{r}^{+}(t)$ and $h_{r}^{+}(t)=\max
\{f_{r}^{-}(0)-f_{r}^{-}(t),f_{r}^{+}(0)-f_{r}^{+}(t)\}=\max
\{-f_{r}^{-}(t),-f_{r}^{+}(t)\}=-f_{r}^{-}(t),$ the limit 
\begin{eqnarray*}
\lim_{t\rightarrow 0^{+}}\frac{h_{r}(t)}{-t} &=&\lim_{t\rightarrow 0^{+}}%
\frac{[h_{r}^{-}(t),h_{r}^{+}(t)]}{-t} \\
&=&\lim_{t\rightarrow 0^{+}}\frac{[-f_{r}^{+}(t),-f_{r}^{-}(t)]}{-t} \\
&=&\lim_{t\rightarrow 0^{+}}\frac{[f_{r}^{-}(t),f_{r}^{+}(t)]}{t} \\
&=&\lim_{t\rightarrow 0^{+}}[\frac{(1+r)t\sin (1/t)}{t},\frac{(3-r)t\sin
(1/t)}{t}] \\
&=&\lim_{t\rightarrow 0^{+}}[(1+r)\sin (1/t),(3-r)\sin (1/t)]
\end{eqnarray*}%
does not exits. Thus, by \cite{Vasavi}, $f$ is not \ $\Delta _{gH}$
differentiable at $0$. To show that the existence of the $\Delta _{gH}$
derivative does not imply the existence of the $\nabla _{gH}$ derivative, we may
consider the same function at $0$ on $\mathbb T=[-3,0]\cup \lbrack 1,2]$.
\end{proof}
\end{theorem}

The following Corollaries are direct consequences of Theorem 4 in \cite{bayeg}.

\begin{corollary}\label{corol2}
Let $s\in \mathbb T_{\kappa }^{\kappa }$ be dense and 
\begin{equation*}
f^{\prime }(s)=\lim_{t\rightarrow s}\frac{f(t)\ominus _{gH}f(s)}{t-s}
\end{equation*}%
exists. Then%
\begin{equation*}
(\diamond _{gH}^{\alpha }f)(s)=(\Delta _{gH}f)(s)=(\nabla
_{gH}f)(s)=f^{\prime }(s),
\end{equation*}%
where $0\leq \alpha \leq 1$.
\end{corollary}

\begin{corollary}\label{corol3}
If $s\in \mathbb T_{\kappa }^{\kappa }$  is a left-scattered and right-dense point and

\begin{equation*}
\left( f^{\prime }(s)\right) _{+}=\lim_{t\rightarrow s^{+}}\frac{%
f(t)\ominus _{gH}f(s)}{t-s}
\end{equation*}%
exists, then 

\begin{enumerate}
\item $(\Delta _{gH}f)(s)=\left( f^{\prime }(s)\right) _{+}$

\item $(\nabla _{gH}f)(s)=\frac{f(\rho (s))\ominus _{gH}f(s)}{\rho (s)-s}$

\item $(\diamond _{gH}^{\alpha }f)(s)=\alpha \left( f^{\prime }(s)\right)
_{+}+(1-\alpha )\frac{f(\rho (s))\ominus _{gH}f(s)}{\rho (s)-s}.$
\end{enumerate}
\end{corollary}

\bigskip

\begin{corollary}\label{corol4}
If $s\in \mathbb T_{\kappa }^{\kappa }$  is a right-scattered and left-dense point and

\begin{equation*}
\left( f^{\prime }(s)\right) _{-}=\lim_{t\rightarrow s^{-}}\frac{%
f(t)\ominus _{gH}f(s)}{t-s}
\end{equation*}%
exists, then 

\begin{enumerate}
\item $(\Delta _{gH}f)(s)=\frac{f(\sigma (s))\ominus _{gH}f(s)}{\sigma (s)-s}%
$

\item $(\nabla _{gH}f)(s)=\left( f^{\prime }(s)\right) _{-}$

\item $(\diamond _{gH}^{\alpha }f)(s)=\alpha \frac{f(\sigma (s))\ominus
_{gH}f(s)}{\sigma (s)-s}+(1-\alpha )\left( f^{\prime }(s)\right) _{-}.$
\end{enumerate}
\end{corollary}

\begin{theorem}\label{thm5}
Let $\mathbb T$ be a time scale, $s\in \mathbb{T}_{\kappa }^{\kappa }$ and $0\leq \alpha \leq
1$. If $f$ is $\diamond _{gH}^{\alpha }$-differentiable at $s,$ then $f$ is continuous at $s$.

\begin{proof}
Assume that $f$ is $\diamond _{gH}^{\alpha }$ differentiable at $s\in \mathbb T_{\kappa
}^{\kappa }$. We will consider the following four cases:

\textbf{Case 1:} Let $s$ be a dense point. Then by Corollary above and \cite{Leelavathi, Vasavi, Fard}
\ we have 
\begin{equation*}
(\diamond _{gH}^{\alpha }f)(s)=(\Delta _{gH}f)(s)=(\nabla _{gH}f)(s)
\end{equation*}

\textbf{Case 2:} Let $s$ be an isolated point. Then $f$ is continuous at $s$.

\textbf{Case 3: } Assume $s$ is right-dense and left-scattered. Thus $\sigma
(s)=s$ and $\rho (s)<s$.

Let $\epsilon \in (0,1)$ and 
\begin{equation*}
\epsilon _{\ast }=\frac{\epsilon \alpha (s-\rho (s))}{(s-\rho
(s)+1)+\left\vert 1-\alpha \right\vert D_{\infty }((f(\rho (s))\ominus
_{gH}f(s)),0)+(s-\rho (s)+1)D_{\infty }((\diamond _{gH}^{\alpha }f)(s),0).}.
\end{equation*}%
Thus $0<\epsilon _{\ast }<1$. Then there is a neighborhood $U_{1}$ of $s$
such that for all $t\in U_{1}$ such that 
\begin{equation*}
D_{\infty }(\alpha \lbrack f(\sigma (s))\ominus _{gH}f(t)]\nu
_{st}+(1-\alpha )[f(\rho (s))\ominus _{gH}f(t)]\mu _{st},(\diamond
_{gH}^{\alpha }f)(s)\mu _{st}\nu _{st})\leq \epsilon _{\ast }\left\vert \mu
_{st}\nu _{st}\right\vert .
\end{equation*}
Since 
\begin{align*}
	D_{\infty} \Big( &\alpha \big[f(\sigma(s)) \ominus_{gH} f(t)\big] \nu_{st} + (1 - \alpha) \big[f(\rho(s)) \ominus_{gH} f(t)\big] \mu_{st}, (\diamond_{gH}^{\alpha} f)(s) \mu_{st} \nu_{st} \Big) \\
	=& D_{\infty} \Big( \alpha \big[f(\sigma(s)) \ominus_{gH} f(t)\big] (\rho(s) - t) + (1 - \alpha) \big[f(\rho(s)) \ominus_{gH} f(t)\big] (\sigma(s) - t), \\
	&\quad (\diamond_{gH}^{\alpha} f)(s) (\sigma(s) - t) (\rho(s) - t) \Big) \\
	=& D_{\infty} \Big( \alpha \big[f(s) \ominus_{gH} f(t)\big] \big((\rho(s) - s) + (s - t)\big) \\
	&\quad + (1 - \alpha) \Big[ \big(f(\rho(s)) \ominus_{gH} f(s)\big) + \big(f(s) \ominus_{gH} f(t)\big) \Big] (s - t), \\
	&\quad (\diamond_{gH}^{\alpha} f)(s) (s - t) \big((\rho(s) - s) + (s - t)\big) \Big) \\
	=& D_{\infty} \Big( \alpha \big[f(s) \ominus_{gH} f(t)\big] (\rho(s) - s) + \alpha \big[f(s) \ominus_{gH} f(t)\big] (s - t) \\
	&\quad + (1 - \alpha) \big(f(\rho(s)) \ominus_{gH} f(s)\big) (s - t) + (1 - \alpha) \big(f(s) \ominus_{gH} f(t)\big) (s - t), \\
	&\quad (\diamond_{gH}^{\alpha} f)(s) (s - t) \big((\rho(s) - s) + (s - t)\big) \Big) \\
	=& D_{\infty} \Big( \big[f(s) \ominus_{gH} f(t)\big] \big( \alpha (\rho(s) - s) + (s - t) \big) \\
	&\quad + (1 - \alpha) \big(f(\rho(s)) \ominus_{gH} f(s)\big) (s - t), \\
	&\quad (\diamond_{gH}^{\alpha} f)(s) (s - t) \big((\rho(s) - s) + (s - t)\big) \Big) \\
	\geq & \left| D_{\infty} \Big( \big[f(s) \ominus_{gH} f(t)\big] \big( \alpha (\rho(s) - s) + (s - t) \big), 0 \Big) \right. \\
	& \quad - \left. D_{\infty} \Big( (1 - \alpha) \big(f(\rho(s)) \ominus_{gH} f(s)\big) (s - t), (\diamond_{gH}^{\alpha} f)(s) (s - t) \big((\rho(s) - s) + (s - t)\big) \Big) \right| \text{.}
\end{align*}
then we have
\begin{align*}
	\left\vert D_{\infty }( [f(s) \ominus _{gH} f(t)] \left( \alpha (\rho (s) - s) + (s - t) \right), 0 ) - D_{\infty } \left( (1 - \alpha ) ( f(\rho (s)) \ominus _{gH} f(s) ) (s - t), \right. \right. \\
	\left. \left. \quad \quad \quad  \left( \diamond _{gH}^{\alpha } f \right)(s) (s - t) ( (\rho (s) - s) + (s - t) ) \right) \right\vert
	 \leq \epsilon_{\ast} \left\vert \mu_{st} \nu_{st} \right\vert \\
	 = \epsilon_{\ast} \left\vert (\sigma(s) - t)(\rho(s) - t) \right\vert \\
	 = \epsilon_{\ast} \left\vert (s - t) \left( (\rho (s) - s) + (s - t) \right) \right\vert .
\end{align*}
Thus,
\begin{align*}
	D_{\infty }&\left([f(s)\ominus _{gH}f(t)](\alpha (\rho (s)-s)+(s-t)),0\right) \\
	&\leq \epsilon _{\ast }\left\vert (s-t)((\rho (s)-s)+(s-t))\right\vert \\
	&\quad + D_{\infty } \left( (1-\alpha )(f(\rho (s))\ominus _{gH}f(s))(s-t), (\diamond _{gH}^{\alpha}f)(s)(s-t)((\rho (s)-s)+(s-t)) \right) \\
	&\leq \epsilon _{\ast }\left\vert (s-t)((\rho (s)-s)+(s-t))\right\vert \\
	&\quad + D_{\infty } \left( (1-\alpha )(f(\rho (s))\ominus _{gH}f(s))(s-t), 0 \right) \\
	&\quad + D_{\infty } \left( (\diamond _{gH}^{\alpha}f)(s)(s-t)((\rho (s)-s)+(s-t)), 0 \right) \\
	&= \epsilon _{\ast }\left\vert (s-t)((\rho (s)-s)+(s-t))\right\vert \\
	&\quad + \left\vert (1-\alpha )(s-t) \right\vert D_{\infty } \left( (f(\rho (s))\ominus _{gH}f(s)), 0 \right) \\
	&\quad + \left\vert (s-t)((\rho (s)-s)+(s-t)) \right\vert D_{\infty } \left( (\diamond _{gH}^{\alpha}f)(s), 0 \right) \\
	&\leq \epsilon _{\ast }(s-\rho (s)+1) \\
	&\quad + \left\vert 1-\alpha \right\vert \epsilon _{\ast } D_{\infty } \left( (f(\rho (s))\ominus _{gH}f(s)), 0 \right) \\
	&\quad + (s-\rho (s)+1) \epsilon _{\ast } D_{\infty } \left( (\diamond _{gH}^{\alpha}f)(s), 0 \right).
\end{align*}

Thus, we obtain 
\begin{align*}
D_{\infty }(f(s)\ominus _{gH}f(t),0)&\leq \epsilon _{\ast }\frac{(s-\rho
(s)+1)+\left\vert 1-\alpha \right\vert D_{\infty }((f(\rho (s))\ominus
_{gH}f(s)),0)}{%
\alpha (s-\rho (s))}\\ &+\epsilon _{\ast }\frac{(s-\rho (s)+1)D_{\infty }((\diamond _{gH}^{\alpha }f)(s),0)}{\alpha (s-\rho (s))}\\
&=\epsilon.
\end{align*}
\end{proof}
\end{theorem}

\begin{theorem}\label{thm6}
Let $\mathbb T$ be a time scale and $0<\alpha <1$. If $f$ is $\diamond _{gH}^{\alpha
}$-differentiable at $s$, then $f$ is both $\Delta _{gH}$ and $\nabla _{gH}$%
-differentiable at $s$.
\end{theorem}

\begin{proof}
Let $\mathbb T$ be a time scale and $0<\alpha <1$. Let $\epsilon >0$ be given, and $%
\epsilon _{\ast }=\epsilon \frac{1-\alpha }{1+\alpha }$.

Assume $f$ is $\diamond _{gH}^{\alpha }$ differentiable at $s\in T_{\kappa
}^{\kappa }.$ We will consider following cases:

\textbf{Case 1: }If $s\in \mathbb T_{\kappa }^{\kappa }$ is dense, then the result
follows from Theorem 4 in \cite{bayeg} and Corollary 2.

\textbf{Case 2: }If $s\in \mathbb T_{\kappa }^{\kappa }$ is isolated, then the
result is obvious.

\textbf{Case 3: }Let $s\in \mathbb T_{\kappa }^{\kappa }$ be right-scattered and
left-dense. Thus $\sigma (s)>s$ and $\rho (s)=s$. Since $f$ is $\diamond
_{gH}^{\alpha }$-differentiable at $s$, $f$ is continuous at $s$. Hence $f$
\ has $\Delta _{gH}$-derivative \cite{Vasavi, Fard}. Then for all $\epsilon _{\ast }>0$
there is a neighbourhood $N_{1}$ of $s$ such that for all $t\in N_{1}:$
\begin{equation*}
D_{\infty }(\alpha \lbrack f(\sigma (s))\ominus _{gH}f(t)]\nu
_{st}+(1-\alpha )[f(t)\ominus _{gH}f(\rho (s))]\mu _{st},(\diamond
_{gH}^{\alpha }f)(s)\mu _{st}\nu _{st})\leq \epsilon _{\ast }\left\vert \mu
_{st}\nu _{st}\right\vert
\end{equation*}
and neighborhood $N_{2}$ of such that $s\in N_{2}$%
\begin{equation*}
D_{\infty }(f(\sigma (s))\ominus _{gH}f(t),(\Delta _{gH}f)(s)\mu _{st})\leq
\epsilon _{\ast }\left\vert \mu _{st}\right\vert .
\end{equation*}
Choose $\gamma $ such that 
\begin{equation*}
(\diamond _{gH}^{\alpha }f)(s)=\alpha (\Delta _{gH}f)(s)+(1-\alpha )\gamma .
\end{equation*}%
Then there exists a neighbourhood $N=N_{1}\cap N_{2}$ of $s$ such that for
all $s\in N:$
\begin{equation*}
\begin{aligned}
\epsilon_{\ast} \left\vert \mu_{st} \nu_{st} \right\vert 
\geq & \, D_{\infty} \Big( \alpha \left[ f(\sigma(s)) \ominus_{gH} f(t) \right] \nu_{st} + (1-\alpha) \left[ f(t) \ominus_{gH} f(\rho(s)) \right] \mu_{st}, \\
& \quad (\diamond_{gH}^{\alpha} f)(s) \mu_{st} \nu_{st} \Big) \\
= & \, D_{\infty} \Big( \alpha \left[ f(\sigma(s)) \ominus_{gH} f(t) \right] \nu_{st} + (1-\alpha) \left[ f(t) \ominus_{gH} f(\rho(s)) \right] \mu_{st}, \\
& \quad \left[ \alpha (\Delta_{gH} f)(s) + (1-\alpha) \gamma \right] \mu_{st} \nu_{st} \Big) \\
\geq & \, \left\vert \alpha \right\vert \left\vert \nu_{st} \right\vert D_{\infty} \Big( f(\sigma(s)) \ominus_{gH} f(t), (\Delta_{gH} f)(s) \mu_{st} \Big) \\
& - (1-\alpha) \left\vert \mu_{st} \right\vert D_{\infty} \Big( f(t) \ominus_{gH} f(\rho(s)), \gamma \nu_{st} \Big) \\
\geq & \, -\alpha \left\vert \nu_{st} \right\vert D_{\infty} \Big( f(\sigma(s)) \ominus_{gH} f(t), (\Delta_{gH} f)(s) \mu_{st} \Big) \\
& + (1-\alpha) \left\vert \mu_{st} \right\vert D_{\infty} \Big( f(t) \ominus_{gH} f(\rho(s)), \gamma \nu_{st} \Big) \\
\geq & \, (1-\alpha) \left\vert \mu_{st} \right\vert D_{\infty} \Big( f(t) \ominus_{gH} f(\rho(s)), \gamma \nu_{st} \Big) \\
& - \alpha \epsilon_{\ast} \left\vert \nu_{st} \right\vert \left\vert \mu_{st} \right\vert.
\end{aligned}
\end{equation*}

\normalsize
This implies%
\begin{eqnarray*}
D_{\infty }(f(t)\ominus _{gH}f(\rho (s)),\gamma \nu _{st}) &\leq &\epsilon
_{\ast }\frac{1+\alpha }{1-\alpha }\left\vert \nu _{st}\right\vert \\
&\leq &\epsilon \left\vert \nu _{st}\right\vert .
\end{eqnarray*}%
Thus $(\nabla _{gH}f)(s)=\gamma $ exists.

\textbf{Case 4:} The case $s$ right-dense, left-scatter is similar.
\end{proof}

\begin{theorem}\label{thm7}
Let $f$,$g:\mathbb T\rightarrow \mathcal F_N(R)$ be $\diamond _{gH}^{\alpha }$-differentiable
at $s\in \mathbb T_{\kappa }^{\kappa }.$ Then we have the followings:

\begin{enumerate}
\item $f+g:\mathbb T\rightarrow \mathcal F_N(R)$ is $\diamond _{gH}^{\alpha }$-differentiable
at $s\in \mathbb T_{\kappa }^{\kappa }$ such that%
\begin{equation*}
(\diamond _{gH}^{\alpha }(f+h))(s)=(\diamond _{gH}^{\alpha }f)(s)+(\diamond
_{gH}^{\alpha }h)(s).
\end{equation*}

\item For any constant $c$, $cf:\mathbb T\rightarrow \mathcal F_N(R)$ is $\diamond
_{gH}^{\alpha }$-differentiable at $s\in \mathbb T_{\kappa }^{\kappa }$ such that%
\begin{equation*}
(\diamond _{gH}^{\alpha }(cf))(s)=c(\diamond _{gH}^{\alpha }f)(s)
\end{equation*}

\item $(\diamond _{gH}^{\alpha }(fh))(s)=(\diamond _{gH}^{\alpha
}f)(s)h(s)+\alpha f(\sigma (s))\Delta _{gH}h(s)+(1-\alpha )f(\rho (s))\nabla
_{gH}h(s)$.
\end{enumerate}
\end{theorem}

\begin{proof}
By \cite{Vasavi, Fard} and \cite{Leelavathi} we have 
\begin{eqnarray*}
\Delta _{gH}(fh)(s) &=&(\Delta _{gH}f)(s)h(s)+f(\sigma (s))(\Delta
_{gH}h)(s), \\
\nabla _{gH}(fh)(s) &=&(\nabla _{gH}f)(s)h(s)+f(\rho (s))(\nabla _{gH}h)(s).
\end{eqnarray*}%
Thus,%
\begin{eqnarray*}
(\diamond _{gH}^{\alpha }(fh))(s) &=&\alpha \Delta _{gH}(fh)(s)+(1-\alpha
)\nabla _{gH}(fh)(s) \\
&=&\alpha ((\Delta _{gH}f)(s)h(s)+f(\sigma (s))(\Delta _{gH}h)(s)) \\
&&+(1-\alpha )((\nabla _{gH}f)(s)h(s)+f(\rho (s))(\nabla _{gH}h)(s)) \\
&=&\alpha (\Delta _{gH}f)(s)h(s)+\alpha f(\sigma (s))(\Delta _{gH}h)(s) \\
&&+(\nabla _{gH}f)(s)h(s)+f(\rho (s))(\nabla _{gH}h)(s) \\
&&-\alpha (\nabla _{gH}f)(s)h(s)-\alpha f(\rho (s))(\nabla _{gH}h)(s) \\
&=&(\alpha (\Delta _{gH}f)(s)+(1-\alpha )(\nabla _{gH}f)(s))h(s)+\alpha
f(\sigma (s))(\Delta _{gH}h)(s) \\
&&+(1-\alpha )(f(\rho (s))(\nabla _{gH}h)(s)).
\end{eqnarray*}
\end{proof}

\section{Fuzzy Diamond-$\protect\alpha $ Integration}

\begin{definition} \label{def15}
A selector of a set valued mapping $G$ from $\mathbb T$ to $K_{C}(%
\mathbb{R}
)$ is a single valued mapping $g:\mathbb T\rightarrow 
\mathbb{R}
$ such that $g(t)\in G(t)$
for all $t\in \mathbb T$.

We will denote the set of all integrable selectors of a set valued function $%
G:\mathbb T\rightarrow K_{C}(%
\mathbb{R}
)$ by $S(G)$. That is,%
\begin{equation*}
S(G)=\left\{ g:\mathbb T\rightarrow 
\mathbb{R}
:g\text{ is integrable and }g(t)\in G(t)\text{ for all }t\in \mathbb T\right\} .
\end{equation*}
\end{definition}

\begin{definition} \label{def16}
The Aumann integral of a fuzzy valued function $f:[a,b]_{\mathbb T}\rightarrow \mathcal F_N(
\mathbb{R})$ is defined levelwise by%
\begin{eqnarray*}
\left[ \int_{a}^{b}f(t)\diamond _{gH}^{\alpha }t\right] _{r} &=&\left\{
\int_{a}^{b}g(t)\diamond _{gH}^{\alpha }t:g\in S( f_{r}(t))\right\} \\
&=&\int_{a}^{b} f_{r}(t)\diamond _{gH}^{\alpha }t
\end{eqnarray*}%
for all $r\in \lbrack 0,1].$
\end{definition}

\begin{definition} \label{def17}
Let $f:[a,b]_{\mathbb T}\rightarrow \mathcal F_N(
\mathbb{R}
)$ be a fuzzy function on time scale $\mathbb T$. We say that $\ f$ is Aumann diamond-alpha integrable on 
$[a,b]_{\mathbb T}$ , which we denote its integral on $[a,b]_{\mathbb T}$ by $%
\int_{a}^{b}f(t)\diamond _{gH}^{\alpha }t$ if there exists a fuzzy number $u\in \mathcal F_N(
\mathbb{R}
)$ such that $[u]_{r}=\int_{a}^{b}f_{r}(t)\diamond _{gH}^{\alpha }t$ for all $r\in \lbrack 0,1]$.
\end{definition}
For simplicity, we will call an Aumann diamond-alpha integrable function as diamond-alpha integrable or $\diamond _{gH}^{\alpha }$-integrable.
\begin{theorem}\label{thm8}
Let $s_{0}$, $s\in \mathbb T$ with $s_{0}<s$, $k\in 
\mathbb{R}$ and $f,g:[s_{0},s]_{\mathbb T}\rightarrow \mathcal F_N(
\mathbb{R})$ be $\diamond _{gH}^{\alpha }$-integrable. Then

\begin{enumerate}
\item $f+g$ is $\diamond _{gH}^{\alpha }$-integrable on $[s_{0},s]_{\mathbb T}$ such
that 
\begin{equation*}
\int_{s_{0}}^{s}\left[ f(t)+g(t)\right] \diamond _{gH}^{\alpha
}t=\int_{s_{0}}^{s}f(t)\diamond _{gH}^{\alpha
}t+\int_{s_{0}}^{s}g(t)\diamond _{gH}^{\alpha }t.
\end{equation*}

\item $kf$ is $\diamond _{gH}^{\alpha }$-integrable on $[s_{0},s]_{\mathbb T}$ such
that 
\begin{equation*}
\int_{s_{0}}^{s}kf(t)\diamond _{gH}^{\alpha }t=k\int_{s_{0}}^{s}f(t)\diamond
_{gH}^{\alpha }t.
\end{equation*}

\item For any $m\in \lbrack s_{0},s]_{\mathbb T}$, we have 
\begin{equation*}
\int_{s_{0}}^{s}f(t)\diamond _{gH}^{\alpha }t=\int_{s_{0}}^{m}f(t)\diamond
_{gH}^{\alpha }t+\int_{m}^{s}f(t)\diamond _{gH}^{\alpha }t.
\end{equation*}
\end{enumerate}
\end{theorem}

\begin{proof}
\begin{enumerate}
\item Since $f$ \ and $g$ are $\diamond _{gH}^{\alpha }$-integrable on $%
[s_{0},s]_{\mathbb T}$, there exists fuzzy numbers $u$ and $v$ in $F(
\mathbb{R}
)$ such that 
\begin{eqnarray*}
 u_{r} &=&\int_{s_{0}}^{s}f_{r}(t)\diamond _{gH}^{\alpha }t, \\
v_{r} &=&\int_{s_{0}}^{s}g_{r}(t)\diamond _{gH}^{\alpha }t
\end{eqnarray*}%
for all $r\in \lbrack 0,1]$. Thus,

\begin{eqnarray*}
\left[ \int_{s_{0}}^{s}\left[ f(t)+g(t)\right] \diamond _{gH}^{\alpha }t%
\right] _{r} &=&\left[ \int_{s_{0}}^{s}\left[ f(t)+g(t)\right] \diamond
_{gH}^{\alpha }t\right] _{r} \\
&=&\int_{s_{0}}^{s}[f(t)+g(t)]_{r}\diamond _{gH}^{\alpha }t \\
&=&\int_{s_{0}}^{s}f_{r}(t)\diamond _{gH}^{\alpha
}t+\int_{s_{0}}^{s}g_{r}(t)\diamond _{gH}^{\alpha }t \\
&=&u_{r}+v_{r} \\
&=&[u+v]_{r}
\end{eqnarray*}%
for all $r\in \lbrack 0,1]$. Since $u+v\in \mathcal F_N(\mathbb{R}
)$, then $f+g$ is also $\diamond _{gH}^{\alpha }$-integrable on $%
[s_{0},s]_{\mathbb T}.$

\item Since $f$ \ is $\diamond _{gH}^{\alpha }$-integrable on $[s_{0},s]_{\mathbb T}$%
, there exists fuzzy numbers $u$ in $\mathcal F_N(
\mathbb{R}
)$ such that 
\begin{equation*}
u_{r}=\int_{s_{0}}^{s}f_{r}(t)\diamond _{gH}^{\alpha }t
\end{equation*}%
for all $r\in \lbrack 0,1]$. Without loss of generality assume that $k>0$.
Thus,

\begin{eqnarray*}
\left[ \int_{s_{0}}^{s}\left[ kf(t)\right] \diamond _{gH}^{\alpha }t\right]
_{r} &=&\int_{s_{0}}^{s}[kf(t)]_{r}\diamond _{gH}^{\alpha }t \\
&=&k\int_{s_{0}}^{s}f_{r}(t)\diamond _{gH}^{\alpha }t \\
&=&ku_{r} \\
&=&[ku]_{r} 
\end{eqnarray*}%
for all $r\in \lbrack 0,1].$ Since $ku\in \mathcal F_N(\mathbb{R})$, then $kf$ is also $\diamond _{gH}^{\alpha }$-integrable on $%
[s_{0},s]_{\mathbb T}.$

\item Let $m\in \lbrack s_{0},s]_{\mathbb T}$. Since $f$ \ is $\diamond
_{gH}^{\alpha }$-integrable on $[s_{0},s]_{\mathbb T}$, then it is $\diamond
_{gH}^{\alpha }$-integrable on $%
[s_{0},m]_{\mathbb T}$ and $[m,s]_{\mathbb T}$. Thus there exists fuzzy numbers $u$ and $v$
in $\mathcal F_N(
\mathbb{R}
)$ such that 
\begin{eqnarray*}
 u_{r} &=&\int_{s_{0}}^{m}f_{r}(t)\diamond _{gH}^{\alpha }t \\
v_{r} &=&\int_{m}^{s}f_{r}(t)\diamond _{gH}^{\alpha }t
\end{eqnarray*}
for all $r\in \lbrack 0,1]$. Thus,

\begin{eqnarray*}
\left[ \int_{s_{0}}^{s}f(t)\diamond _{gH}^{\alpha }t\right] _{r}
&=&\int_{s_{0}}^{s}f_{r}(t)\diamond _{gH}^{\alpha }t \\
&=&\int_{s_{0}}^{m}f_{r}(t)\diamond _{gH}^{\alpha
}t+\int_{m}^{s}f_{r}(t)\diamond _{gH}^{\alpha }t \\
&=&\left[ \int_{s_{0}}^{m}f(t)\diamond _{gH}^{\alpha }t\right] _{r}+\left[
\int_{m}^{s}f(t)\diamond _{gH}^{\alpha }t\right] _{r} \\
&=&\left[ \int_{s_{0}}^{m}f(t)\diamond _{gH}^{\alpha
}t+\int_{m}^{s}f(t)\diamond _{gH}^{\alpha }t\right] _{r}
\end{eqnarray*}

for all $r\in \lbrack 0,1]$. Therefore, for any $m\in \lbrack s_{0},s]_{\mathbb T}$,
we have 
\begin{equation*}
\int_{s_{0}}^{s}f(t)\diamond _{gH}^{\alpha }t=\int_{s_{0}}^{m}f(t)\diamond
_{gH}^{\alpha }t+\int_{m}^{s}f(t)\diamond _{gH}^{\alpha }t.
\end{equation*}
\end{enumerate}
\end{proof}

\begin{theorem}\label{thm9}
If $f,h:[a,b]_{\mathbb T}\rightarrow \mathcal F_N(
\mathbb{R})$ are $\diamond _{gH}^{\alpha }$-integrable, then $D_{\infty }(f,h)$ is
also $\diamond
_{gH}^{\alpha }$-integrable on $[a,b]_{\mathbb T}$ and 
\[
D_{\infty }\left( \int_{a}^{b}f(t)\diamond _{gH}^{\alpha
}t,\int_{a}^{b}h(t)\diamond _{gH}^{\alpha }t\right) \leq
\int_{a}^{b}D_{\infty }\left( f(t),h(t\right) )\diamond ^{\alpha }t
\]
\end{theorem}

\begin{proof}
Since%
\begin{align*}
D\left(\int_{a}^{b}f_{r}(t)\diamond_{gH}^{\alpha}t, 
\int_{a}^{b}h_{r}(t)\diamond_{gH}^{\alpha}t\right) 
&= D\left( \int_{a}^{b}\left[ f_{r}^{-}(t), f_{r}^{+}(t) \right] 
\diamond_{gH}^{\alpha}t, \int_{a}^{b}\left[ g_{r}^{-}(t), g_{r}^{+}(t) \right] 
\diamond_{gH}^{\alpha}t \right) \\
&= D\left( \left[ \int_{a}^{b}f_{r}^{-}(t)\diamond^{\alpha}t, 
\int_{a}^{b}f_{r}^{+}(t)\diamond^{\alpha}t \right], \right. \\
&\quad \left. \left[ \int_{a}^{b}g_{r}^{-}(t)\diamond^{\alpha}t, 
\int_{a}^{b}g_{r}^{+}(t)\diamond^{\alpha}t \right] \right) \\
&= \max\left\{ \left\vert \int_{a}^{b}f_{r}^{-}(t)\diamond^{\alpha}t 
- \int_{a}^{b}g_{r}^{-}(t)\diamond^{\alpha}t \right\vert, \right. \\
&\quad \left. \left\vert \int_{a}^{b}f_{r}^{+}(t)\diamond^{\alpha}t 
- \int_{a}^{b}g_{r}^{+}(t)\diamond^{\alpha}t \right\vert \right\} \\
&= \max\left\{ \left\vert \int_{a}^{b}\left( f_{r}^{-}(t) - g_{r}^{-}(t) 
\right) \diamond^{\alpha}t \right\vert, \right. \\
&\quad \left. \left\vert \int_{a}^{b}\left( f_{r}^{+}(t) - g_{r}^{+}(t) 
\right) \diamond^{\alpha}t \right\vert \right\} \\
&\leq \max\left\{ \int_{a}^{b}\left\vert f_{r}^{-}(t) - g_{r}^{-}(t) 
\right\vert \diamond^{\alpha}t, \right. \\
&\quad \left. \int_{a}^{b}\left\vert f_{r}^{+}(t) - g_{r}^{+}(t) 
\right\vert \diamond^{\alpha}t \right\} \\
&\leq \int_{a}^{b} \max\left\{ \left\vert f_{r}^{-}(t) - g_{r}^{-}(t) 
\right\vert, \left\vert f_{r}^{+}(t) - g_{r}^{+}(t) \right\vert \right\} 
\diamond^{\alpha}t \\
&= \int_{a}^{b}D\left( f_{r}(t), h_{r}(t) \right) \diamond^{\alpha}t.
\end{align*}

then we have%
\[
D_{\infty }\left( \int_{a}^{b}f(t)\diamond _{gH}^{\alpha
}t,\int_{a}^{b}h(t)\diamond _{gH}^{\alpha }t\right) \leq
\int_{a}^{b}D_{\infty }\left( f(t),h(t\right) )\diamond ^{\alpha }t.
\]
\end{proof}
\newpage
\begin{lemma}\label{lem2}
Let $f:[a,b]_{\mathbb T}\rightarrow \mathcal F_N(
\mathbb{R}
)$ be a fuzzy function on time scale $\mathbb T$. If $f$ is regulated at $t$, then
the end-points of the r-level sets $f_{r}^{-}$ and $f_{r}^{+}$ are regulated
at $t$ for any $r\in \lbrack 0,1]$.
\end{lemma}

For the convenience of the following result, we define the following sets:
\begin{eqnarray*}
S_{1} &=&\left\{ t\in \mathbb T:t\text{ is left-dense and right-scattered}\right\} ,
\\
S_{2} &=&\left\{ t\in \mathbb T:t\text{ is left-scattered and right-dense}\right\} ,
\\
S_{3} &=&\left\{ t\in \mathbb T:t\text{ is left-scattered and right-scattered}%
\right\} , \\
S_{4} &=&\left\{ t\in \mathbb T:t\text{ is left-dense and right-dense}\right\} .
\end{eqnarray*}

\begin{theorem}\label{thm10}
Let $h:[a,b]_{\mathbb T}\rightarrow \mathcal F_N(
\mathbb{R})$ be a fuzzy function on a time scale $\mathbb T$. If $h$ is regulated in $\mathbb T$, then $%
\int_{t_{0}}^{t}h(s)\diamond _{gH}^{\alpha }s$ is $\diamond _{gH}^{\alpha } 
$-differentiable on $\mathbb T$, and%
\begin{equation*}
\begin{split}
\diamond_{gH}^{\alpha} \left( \int_{t_{0}}^{t} h(s) \diamond_{gH}^{\alpha} s \right) 
&= (1 - 2\alpha + 2\alpha^2) h(t) \\
&\quad + \alpha (\alpha - 1) \left\{
\begin{array}{cc}
\lim_{s\rightarrow t^{-}}h(s)+h(\sigma (t)), & t\in S_{1}, \\ 
\lim_{s\rightarrow t^{+}}h(s)+h(\rho (t)), & t\in S_{2}, \\ 
h(\sigma (t))+h(\rho (t)), & t\in S_{3}, \\ 
\lim_{s\rightarrow t^{-}}h(s)+\lim_{s\rightarrow t^{+}}h(s), & t\in S_{4}.%
\end{array}
\right.
\end{split}
\end{equation*}

\end{theorem}

Let us define%
\begin{equation*}
f(t)=\int_{t_{0}}^{t}h(s)\diamond _{gH}^{\alpha }s\text{.}
\end{equation*}

$f$ is $\diamond _{gH}^{\alpha }$-differentiable and 
\begin{eqnarray*}
\left( \diamond _{gH}^{\alpha }f\right) (t) &=&\diamond _{gH}^{\alpha
}\left( \int_{t_{0}}^{t}h(s)\diamond _{gH}^{\alpha }s\right) . \\
\left( \diamond _{gH}^{\alpha }f\right) _{r}(t) &=&\left( \diamond
_{gH}^{\alpha }\left( \int_{t_{0}}^{t}h(s)\diamond _{gH}^{\alpha }s\right)
\right) _{r}
\end{eqnarray*}%
We have four cases:

\textbf{i)} 
\begin{equation*}
\left[ (\diamond ^{\alpha }f_{r}^{-})(t),\left( \diamond ^{\alpha
}f_{r}^{+}\right) (t)\right] =\left[ \diamond ^{\alpha }\left(
\int_{t_{0}}^{t}h_{r}^{-}(s)\diamond ^{\alpha }s\right) ,\diamond ^{\alpha
}\left( \int_{t_{0}}^{t}h_{r}^{+}(s)\diamond ^{\alpha }s\right) \right]
\end{equation*}

\textbf{ii) }
\begin{equation*}
\left[ (\diamond ^{\alpha }f_{r}^{+})(t),\left( \diamond ^{\alpha
}f_{r}^{-}\right) (t)\right] =\left[ \diamond ^{\alpha }\left(
\int_{t_{0}}^{t}h_{r}^{-}(s)\diamond ^{\alpha }s\right) ,\diamond ^{\alpha
}\left( \int_{t_{0}}^{t}h_{r}^{+}(s)\diamond ^{\alpha }s\right) \right]
\end{equation*}

\textbf{iii) }
\begin{equation*}
\left[ (\diamond ^{\alpha }f_{r}^{-})(t),\left( \diamond ^{\alpha
}f_{r}^{+}\right) (t)\right] =\left[ \diamond ^{\alpha }\left(
\int_{t_{0}}^{t}h_{r}^{+}(s)\diamond ^{\alpha }s\right) ,\diamond ^{\alpha
}\left( \int_{t_{0}}^{t}h_{r}^{-}(s)\diamond ^{\alpha }s\right) \right]
\end{equation*}

\textbf{iv) }
\begin{equation*}
\left[ (\diamond ^{\alpha }f_{r}^{+})(t),\left( \diamond ^{\alpha
}f_{r}^{-}\right) (t)\right] =\left[ \diamond ^{\alpha }\left(
\int_{t_{0}}^{t}h_{r}^{+}(s)\diamond ^{\alpha }s\right) ,\diamond ^{\alpha
}\left( \int_{t_{0}}^{t}h_{r}^{-}(s)\diamond ^{\alpha }s\right) \right]
\end{equation*}

\begin{proof}
Let us prove for the first one, the others can be shown similarly.

By \cite{Sheng}, since we have 
\begin{equation*}
\begin{split}
(\diamond^{\alpha} f_{r}^{-})(t) &= \diamond^{\alpha} \left( \int_{t_{0}}^{t} h_{r}^{-}(s) \diamond^{\alpha} s \right) \\
&= (1 - 2\alpha + 2\alpha^2) h_{r}^{-}(t) \\
&\quad + \alpha (\alpha - 1) \left\{
\begin{array}{l}
\lim_{s \to t^{-}} h_{r}^{-}(s) + h_{r}^{-}(\sigma(t)), \quad t \in S_1, \\
\lim_{s \to t^{+}} h_{r}^{-}(s) + h_{r}^{-}(\rho(t)), \quad t \in S_2, \\
h_{r}^{-}(\sigma(t)) + h_{r}^{-}(\rho(t)), \quad t \in S_3, \\
\lim_{s \to t^{-}} h_{r}^{-}(s) + \lim_{s \to t^{+}} h_{r}^{-}(s), \quad t \in S_4
\end{array}
\right.
\end{split}
\end{equation*}

and 
\begin{equation*}
\begin{split}
(\diamond^{\alpha} f_{r}^{+})(t) &= \diamond^{\alpha} \left( \int_{t_{0}}^{t} h_{r}^{+}(s) \diamond^{\alpha} s \right) \\
&= (1 - 2\alpha + 2\alpha^2) h_{r}^{+}(t) \\
&\quad + \alpha (\alpha - 1) \left\{
\begin{array}{l}
\lim_{s \to t^{+}} h_{r}^{+}(s) + h_{r}^{+}(\sigma(t)), \quad t \in S_1, \\
\lim_{s \to t^{+}} h_{r}^{+}(s) + h_{r}^{+}(\rho(t)), \quad t \in S_2, \\
h_{r}^{+}(\sigma(t)) + h_{r}^{+}(\rho(t)), \quad t \in S_3, \\
\lim_{s \to t^{+}} h_{r}^{+}(s) + \lim_{s \to t^{+}} h_{r}^{+}(s), \quad t \in S_4
\end{array}
\right.
\end{split}
\end{equation*}

it follows that
\begin{eqnarray*}
&&\left[ (\diamond ^{\alpha }f_{r}^{-})(t),(\diamond ^{\alpha }f_{r}^{+})(t)%
\right] =\left[ \diamond ^{\alpha }\left(
\int_{t_{0}}^{t}h_{r}^{-}(s)\diamond ^{\alpha }s\right) ,\diamond ^{\alpha
}\left( \int_{t_{0}}^{t}h_{r}^{+}(s)\diamond ^{\alpha }s\right) \right] \\
&=&(1-2\alpha +2\alpha ^{2})\left[ h_{r}^{-}(t),h_{r}^{+}(t)\right] \\ &+& \alpha
(\alpha -1)\left\{ 
\begin{array}{cc}
\left[ \lim_{s\rightarrow t^{-}}h_{r}^{-}(s)+h_{r}^{-}(\sigma
(t)),\lim_{s\rightarrow t^{-}}h_{r}^{+}(s)+h_{r}^{+}(\sigma (t))\right] & 
t\in S_{1}, \\ 
\left[ \lim_{s\rightarrow t^{+}}h_{r}^{-}(s)+h_{r}^{-}(\rho
(t)),\lim_{s\rightarrow t^{+}}h_{r}^{+}(s)+h_{r}^{+}(\rho (t))\right] & t\in
S_{2}, \\ 
\left[ h_{r}^{-}(\sigma (t))+h_{r}^{-}(\rho (t)),h_{r}^{+}(\sigma
(t))+h_{r}^{+}(\rho (t))\right] & t\in S_{3}, \\ 
\left[ \lim_{s\rightarrow t^{-}}h_{r}^{-}(s)+\lim_{s\rightarrow
t^{+}}h_{r}^{-}(s),\lim_{s\rightarrow t^{-}}h_{r}^{+}(s)+\lim_{s\rightarrow
t^{+}}h_{r}^{+}(s)\right] & t\in S_{4}.%
\end{array}%
\right.
\end{eqnarray*}
for any $r\in \lbrack 0,1].$
Hence, we obtain 
\begin{equation*}
\begin{split}
\diamond_{gH}^{\alpha} \left( \int_{t_{0}}^{t} h(s) \diamond_{gH}^{\alpha} s \right) 
&= (1 - 2\alpha + 2\alpha^2) h(t) \\
&\quad + \alpha (\alpha - 1) \left\{
\begin{array}{cc}
\lim_{s\rightarrow t^{-}}h(s)+h(\sigma (t)), & t\in S_{1}, \\ 
\lim_{s\rightarrow t^{+}}h(s)+h(\rho (t)), & t\in S_{2}, \\ 
h(\sigma (t))+h(\rho (t)), & t\in S_{3}, \\ 
\lim_{s\rightarrow t^{-}}h(s)+\lim_{s\rightarrow t^{+}}h(s), & t\in S_{4}.%
\end{array}
\right.
\end{split}
\end{equation*}
\end{proof}

\begin{remark}\label{rem4}
In the previous theorem, if $h$ is continuous, then we have 
\begin{equation*}
\diamond _{gH}^{\alpha }\left( \int_{t_{0}}^{t}h(s)\diamond _{gH}^{\alpha
}s\right) =(1-2\alpha +2\alpha ^{2})h(t)+\alpha (\alpha -1)h(\sigma
(t))+h(\rho (t))
\end{equation*}%
for any $t\in S_{1}$, $S_{2}$, $S_{3}$, and $S_{4}.$
\end{remark}

\begin{corollary}\label{corol5}
Let $t\in \mathbb T_{\kappa }^{\kappa }$ and $f,g:\mathbb T\rightarrow \mathcal F_N(\mathbb{R})$, then

\begin{enumerate}
\item $\int_{t}^{\sigma (t)}f(s)\diamond _{gH}^{\alpha }s=\mu (t)(\alpha
f(t)+(1-\alpha )f(\sigma (t))).$

\item $\int_{\rho (t)}^{t}f(s)\diamond _{gH}^{\alpha }s=\nu (t)(\alpha
f(\rho (t))+(1-\alpha )f(t)).$
\end{enumerate}
\end{corollary}

\begin{proof}
Let $\int_{t}^{\sigma (t)}f(s)\diamond _{gH}^{\alpha }s=u(t)\in \mathcal F_N(
\mathbb{R}
)$. For any $r\in \lbrack 0,1]$, we have
\begin{eqnarray*}
u_{r}(t) &=&\left[ \int_{t}^{\sigma (t)}f(s)\diamond
_{gH}^{\alpha }s\right] _{r} \\
&=&\int_{t}^{\sigma (t)} f_{r}(s)\diamond _{gH}^{\alpha }s \\
&=&\int_{t}^{\sigma (t)}\left[ f_{r}^{-}(s),f_{r}^{+}(s)\right] \diamond
_{gH}^{\alpha }s \\
&=&\left[ \int_{t}^{\sigma (t)}f_{r}^{-}(s)\diamond ^{\alpha
},\int_{t}^{\sigma (t)}f_{r}^{+}(s)\diamond ^{\alpha }\right] \\
&=&\left[ \mu (t)(\alpha f_{r}^{-}(t)+(1-\alpha )f_{r}^{-}(\sigma (t))),\mu
(t)(\alpha f_{r}^{+}(t)+(1-\alpha )f_{r}^{+}(\sigma (t)))\right] \\
&=&\mu (t)\alpha \left[ f_{r}^{-}(t),f_{r}^{+}(t)\right] +\mu (t)(1-\alpha )%
\left[ f_{r}^{-}(t),f_{r}^{+}(t)\right] \\
&=&\mu (t)\alpha f_{r}(t)+\mu (t)(1-\alpha )f_{r}(t).
\end{eqnarray*}%
Hence, we obtain 
\begin{equation*}
\int_{t}^{\sigma (t)}f(s)\diamond _{gH}^{\alpha }s=\mu (t)(\alpha
f(t)+(1-\alpha )f(\sigma (t))).
\end{equation*}%
The proof of (2) can be done in a similar way.
\end{proof}
Now we will give the Fundamental Theorem of Calculus for diamond-alpha integrals of fuzzy functions on a quasi-regular time scale.
\begin{theorem}\label{thm11}
Suppose that $\mathbb T$ is a quasi-regular time scale with 
\begin{equation*}
(\nabla \sigma )(t)=\left\{ 
\begin{array}{cc}
m_{1}, & t<0, \\ 
m_{2}, & t\geq 0.%
\end{array}%
\right. 
\end{equation*}%
$F:\mathbb T\rightarrow  \mathcal F_N(\mathbb{R})$ be $\diamond _{gH}^{\alpha }$-differentiable on an interval $[a,b]_{\mathbb T}$.
Then

\begin{enumerate}
\item If $b<0$, then 
\begin{equation*}
\int\limits_{a}^{b}(\diamond _{gH}^{\alpha }F)(t)\diamond _{gH}^{\alpha
}t=\left. (\alpha ^{2}+(1-\alpha )^{2})F(t)\right\vert _{a}^{b}+\left.
\alpha (1-\alpha )[\frac{F(\sigma (t))}{m_{1}}+m_{1}F(\rho (t))]\right\vert
_{a}^{b}.
\end{equation*}

\item If $a\geq 0$, then 
\begin{equation*}
\int\limits_{a}^{b}(\diamond _{gH}^{\alpha }F)(t)\diamond _{gH}^{\alpha
}t=\left. (\alpha ^{2}+(1-\alpha )^{2})F(t)\right\vert _{a}^{b}+\left.
\alpha (1-\alpha )[\frac{F(\sigma (t))}{m_{2}}+m_{2}F(\rho (t))]\right\vert
_{a}^{b}.
\end{equation*}

\item If $c\in \mathbb T$ is the smallest non-negative point such that $a<c<b$, then
we have 
\begin{eqnarray*}
\int\limits_{a}^{b}(\diamond _{gH}^{\alpha }F)(t)\diamond _{gH}^{\alpha }t
&=&\left. (\alpha ^{2}+(1-\alpha )^{2})F(t)\right\vert _{a}^{c}+\left.
\alpha (1-\alpha )[\frac{F(\sigma (t))}{m_{1}}+m_{1}F(\rho (t))]\right\vert
_{a}^{c} \\
&&+\left. (\alpha ^{2}+(1-\alpha )^{2})F(t)\right\vert _{c}^{b}+\left.
\alpha (1-\alpha )[\frac{F(\sigma (t))}{m_{2}}+m_{2}F(\rho (t))]\right\vert
_{c}^{b}.
\end{eqnarray*}
\end{enumerate}
\end{theorem}

\begin{proof}
Without loss of generality we will assume $F$ is $\diamond _{gH}^{\alpha }$%
-differentiable function on $[a,b]_{T}$ such that $\left[ (\diamond
_{gH}^{\alpha }F)(t)\right] _{r}=\left[ (\diamond ^{\alpha
}F_{r}^{-})(t),(\diamond ^{\alpha }F_{r}^{+})(t)\right] $ for any fixed $%
r\in \lbrack 0,1]$.

\begin{enumerate}
\item Let $b<0$ and $\mathbb T$ be a regular time scale with $(\nabla \sigma )(t)$ $%
=m_{1}\in \lbrack 0,1)$. Then, we have 
\begin{eqnarray*}
\left[ \int\limits_{a}^{b} (\diamond _{gH}^{\alpha }F)(t) \diamond _{gH}^{\alpha } t \right] _{r} 
&=& \int\limits_{a}^{b} \left[ (\diamond _{gH}^{\alpha }F)(t) \right] _{r} 
\diamond _{gH}^{\alpha } t \\
&=& \alpha \int\limits_{a}^{b} \left[ (\diamond _{gH}^{\alpha }F)(t) \right] _{r} \Delta _{gH} t \\
&& + (1-\alpha) \int\limits_{a}^{b} \left[ (\diamond _{gH}^{\alpha }F)(t) \right] _{r} \nabla _{gH} t \\
&=& \alpha \int\limits_{a}^{b} 
\left[ \alpha (\Delta _{gH}F)(t) + (1-\alpha)(\nabla _{gH}F)(t) \right] _{r} \Delta _{gH} t \\
&& + (1-\alpha) \int\limits_{a}^{b} 
\left[ \alpha (\Delta _{gH}F)(t) + (1-\alpha)(\nabla _{gH}F)(t) \right] _{r} \nabla _{gH} t \\
&=& \alpha^2 \int\limits_{a}^{b} \left[ (\Delta _{gH}F)(t) \right] _{r} \Delta _{gH} t \\
&& + (1-\alpha)^2 \int\limits_{a}^{b} \left[ (\nabla _{gH}F)(t) \right] _{r} \nabla _{gH} t  + H(\alpha) \\
&=& (\alpha^2 + (1-\alpha)^2) (F_{r}(b) \ominus _{gH} F_{r}(a)) + H(\alpha),
\end{eqnarray*}
where
\begin{equation*}
H(\alpha) = \alpha (1-\alpha) \left( 
    \int\limits_{a}^{b} \left[ (\Delta _{gH}F)(t) \right] _{r} \nabla _{gH} t 
    + \int\limits_{a}^{b} \left[ (\nabla _{gH}F)(t) \right] _{r} \Delta _{gH} t 
\right).
\end{equation*}
Now, let us consider the function $H(\alpha )$ for $0\leq \alpha \leq 1$. By
the Theorem 3, we obtain 
\begin{eqnarray*}
\int\limits_{a}^{b} \left[ (\Delta _{gH}F)(t) \right] _{r} \nabla _{gH} t 
&=& \int\limits_{a}^{b} \left[ (\Delta F_{r}^{-})(t), (\Delta F_{r}^{+})(t) \right] \nabla _{gH} t \\
&=& \left[ 
    \int\limits_{a}^{b} (\Delta F_{r}^{-})(t) \nabla t, 
    \int\limits_{a}^{b} (\Delta F_{r}^{+})(t) \nabla t 
\right] \\
&=& \left[ 
    \int\limits_{a}^{b} (\Delta F_{r}^{-})(\sigma (t)) \nabla t, 
    \int\limits_{a}^{b} (\Delta F_{r}^{+})(\sigma (t)) \nabla t 
\right] \\
&=& \left[ 
    \int\limits_{a}^{b} (\Delta F_{r}^{-})(\sigma (t)) (\nabla \sigma)(t) (\Delta \rho)(t) \nabla t, \right. \\
&& \quad \left. 
    \int\limits_{a}^{b} (\Delta F_{r}^{+})(\sigma (t)) (\nabla \sigma)(t) (\Delta \rho)(t) \nabla t 
\right] \\
&=& \frac{1}{m_{1}} \left[ 
    \int\limits_{a}^{b} \nabla (F_{r}^{-} \circ \sigma)(t) \nabla t, 
    \int\limits_{a}^{b} \nabla (F_{r}^{+} \circ \sigma)(t) \nabla t 
\right] \\
&=& \frac{1}{m_{1}} \left[ 
    \int\limits_{a}^{b} \nabla (F_{r}^{-} \circ \sigma)(t) \nabla t, 
    \int\limits_{a}^{b} \nabla (F_{r}^{+} \circ \sigma)(t) \nabla t 
\right] \\
&=& \frac{1}{m_{1}} \left[ 
    F_{r}^{-}(\sigma (b)) - F_{r}^{-}(\sigma (a)), 
    F_{r}^{+}(\sigma (b)) - F_{r}^{+}(\sigma (a)) 
\right] \\
&=& \frac{F_{r}(\sigma (b)) \ominus _{gH} F_{r}(\sigma (a))}{m_{1}}
\end{eqnarray*}

and similarly, one can obtain 
\begin{equation*}
\int\limits_{a}^{b}\left[ (\nabla _{gH}F)(t)\right] _{r}\Delta
_{gH}t=m_{1}(F_{r}(\rho (b))\ominus _{gH}F_{r}(\rho (a))).
\end{equation*}%
Therefore, for any fixed $r\in \lbrack 0,1]$ we obtain%
\begin{equation*}
\left[ \int\limits_{a}^{b}(\diamond _{gH}^{\alpha }F)(t)\diamond
_{gH}^{\alpha }t\right] _{r}=(\alpha ^{2}+(1-\alpha )^{2})(F_{r}(b)\ominus
_{gH}F_{r}(a))+H(\alpha ),
\end{equation*}%
where 
\begin{equation*}
H(\alpha )=\alpha (1-\alpha )\left( \frac{F_{r}(\sigma (b))\ominus
_{gH}F_{r}(\sigma (a))}{m_{1}}+m_{1}(F_{r}(\rho (b))\ominus _{gH}F_{r}(\rho
(a)))\right) .
\end{equation*}
Thus, 
\begin{equation*}
\int\limits_{a}^{b}(\diamond _{gH}^{\alpha }F)(t)\diamond _{gH}^{\alpha
}t=\left. (\alpha ^{2}+(1-\alpha )^{2})F(t)\right\vert _{a}^{b}+\left.
\alpha (1-\alpha )[\frac{F(\sigma (t))}{m_{1}}+m_{1}F(\rho (t))]\right\vert
_{a}^{b}.
\end{equation*}

\item Let $a\geq 0$, $\mathbb T$ be a regular time scale with $(\nabla \sigma )(t)$ $%
=m_{2}\geq 1$. Similarly, one can show that 
\begin{equation*}
\left[ \int\limits_{a}^{b}(\diamond _{gH}^{\alpha }F)(t)\diamond
_{gH}^{\alpha }t\right] _{r}=(\alpha ^{2}+(1-\alpha )^{2})(F_{r}(b)\ominus
_{gH}F_{r}(a))+H(\alpha ),
\end{equation*}%
where 
\begin{equation*}
H(\alpha )=\alpha (1-\alpha )\left( \frac{F_{r}(\sigma (b))\ominus
_{gH}F_{r}(\sigma (a))}{m_{2}}+m_{2}(F_{r}(\rho (b))\ominus _{gH}F_{r}(\rho
(a)))\right)
\end{equation*}
for any fixed $r\in \lbrack 0,1]$.

\item Let $c\in \mathbb T$ be the smallest non-negative point such that $a<c<b$, $\mathbb T$
be a regular time scale with $(\nabla \sigma )(t)$ $=m_{2}\in \lbrack 0,1)$
for $a<t<c$ \ and $(\nabla \sigma )(t)$ $=m_{1}\geq 1$ for $c<t<b$. By using
the results in (1) and (2) one can obtain that
\begin{equation*}
\begin{aligned}
\int\limits_{a}^{b} (\diamond _{gH}^{\alpha }F)(t) \diamond _{gH}^{\alpha } t 
=& \left. (\alpha^2 + (1-\alpha)^2) F(t) \right|_{a}^{b} \\
&+ \left. \alpha (1-\alpha) 
\left[ 
    \frac{F(\sigma (t))}{m_{2}} + m_{2} F(\rho (t)) 
\right] \right|_{a}^{c} \\
&+ \left. \alpha (1-\alpha) 
\left[ 
    \frac{F(\sigma (t))}{m_{1}} + m_{1} F(\rho (t)) 
\right] \right|_{c}^{b}.
\end{aligned}
\end{equation*}

\end{enumerate}
\end{proof}

\bigskip

\begin{corollary}\label{corol6}
Suppose that $\mathbb T$ is a homogeneous time scale with $\mu (t)=\nu (t)=k\geq 0$.
\ $F:\mathbb T\rightarrow \mathcal F_N (
\mathbb{R})$ be $\diamond _{gH}^{\alpha }$-differentiable on an interval $[a,b]_\mathbb T$. Then%
\begin{equation*}
\int\limits_{a}^{b}(\diamond _{gH}^{\alpha }F)(t)\diamond _{gH}^{\alpha
}t=\left. (\alpha ^{2}+(1-\alpha )^{2})F(t)\right\vert _{a}^{b}+\left.
\alpha (1-\alpha )(F(t+k)+F(t-k))\right\vert _{c}^{b}.
\end{equation*}
\end{corollary}

The following theorem is about integration by parts for $\diamond
_{gH}^{\alpha }$-integrals of fuzzy functions on a quasi-regular time scale.

\begin{theorem}\label{thm12}
Let $f$ and $g$ be $\diamond _{gH}^{\alpha }$-differentiable functions at $%
t\in \mathbb T_{\kappa }^{\kappa }$ such that $(\diamond _{gH}^{\alpha }f)g$, $%
f(\sigma (t))(\Delta _{gH}g)$ and $f(\rho (t))(\nabla _{gH}g)$ are $\diamond
_{gH}^{\alpha }$-integrable. Then,

\begin{enumerate}
\item If $a\geq 0$, then we have%
\begin{equation*}
\begin{aligned}
\int\limits_{a}^{b} (\diamond _{gH}^{\alpha }f)(t) g(t) \diamond _{gH}^{\alpha } t 
=& \left. (\alpha^2 + (1-\alpha)^2) (fg)(t) \right|_{a}^{b} \\
&+ \left. \alpha (1-\alpha) 
\left[ 
    (fg)(\sigma(t)) + (fg)(\rho(t)) 
\right] \right|_{a}^{b} \\
&\ominus _{gH} 
\alpha \int\limits_{a}^{b} f(\sigma(t)) (\Delta _{gH}g)(t) \diamond _{gH}^{\alpha } t \\
&\ominus _{gH} 
(1-\alpha) \int\limits_{a}^{b} f(\rho(t)) (\nabla _{gH}g)(t) \diamond _{gH}^{\alpha } t.
\end{aligned}
\end{equation*}

\item If $b<0$, then we have%
\begin{equation*}
\begin{aligned}
\int\limits_{a}^{b} (\diamond _{gH}^{\alpha }f)(t) g(t) \diamond _{gH}^{\alpha } t 
=& \left. (\alpha^2 + (1-\alpha)^2) (fg)(t) \right|_{a}^{b} \\
&+ \left. \alpha (1-\alpha) 
\left[ 
    (fg)(\sigma (t)) + (fg)(\rho (t)) 
\right] \right|_{a}^{b} \\
&\ominus _{gH} \alpha \int\limits_{a}^{b} f(\sigma (t)) (\Delta _{gH}g)(t) 
\diamond _{gH}^{\alpha } t \\
&\ominus _{gH} (1-\alpha) \int\limits_{a}^{b} f(\rho (t)) (\nabla _{gH}g)(t) 
\diamond _{gH}^{\alpha } t.
\end{aligned}
\end{equation*}

\item If $c\in \mathbb T$ is the smallest non-negative point such that $a<c<b$, then
we have%
\begin{equation*}
\begin{aligned}
\int\limits_{a}^{b} (\diamond _{gH}^{\alpha }f)(t) g(t) \diamond _{gH}^{\alpha } t 
=& \left. (\alpha^2 + (1-\alpha)^2) (fg)(t) \right|_{a}^{b} \\
&+ \left. \alpha (1-\alpha) 
\left[ 
    (fg)(\sigma(t)) + (fg)(\rho(t)) 
\right] \right|_{a}^{b} \\
&\ominus _{gH} 
\alpha \int\limits_{a}^{b} 
f(\sigma(t)) (\Delta _{gH}g)(t) \diamond _{gH}^{\alpha } t \\
&\ominus _{gH} 
(1-\alpha) \int\limits_{a}^{b} 
f(\rho(t)) (\nabla _{gH}g)(t) \diamond _{gH}^{\alpha } t.
\end{aligned}
\end{equation*}

\end{enumerate}
\end{theorem}

\begin{proof}
By (3) in Theorem \ref{thm7} for $\diamond _{gH}^{\alpha }$-differentiable
functions, we have that 
\begin{equation*}
(\diamond _{gH}^{\alpha }f)(s)h(s)=(\diamond _{gH}^{\alpha }(fh))(s)\ominus
_{gH}\alpha f(\sigma (s))\Delta _{gH}h(s)\ominus _{gH}(1-\alpha )f(\rho
(s))\nabla _{gH}h(s).
\end{equation*}%
If we take the $\diamond _{gH}^{\alpha }$-integral of this equation from $a$
to $b$ and use Theorem \ref{thm11} we obtain the results.
\end{proof}

\section{Conclusions}

In this paper, we have extended the framework of fuzzy calculus on time scales by giving more results on the generalized Hukuhara diamond-alpha derivative. We have introduced generalized Hukuhara diamond-alpha integral for fuzzy number-valued functions. Through the definition and analysis of the generalized Hukuhara diamond-alpha integral, we established key theoretical properties, including its connection to the Fundamental Theorem of Calculus and the method of Integration by Parts.

Our findings demonstrate that the generalized Hukuhara diamond-alpha derivative ensures continuity at points of differentiability and implies both delta and nabla differentiability, offering a more comprehensive framework for analyzing dynamic systems with uncertainty. Similarly, the integration methods we proposed pave the way for addressing complex dynamic models that require handling fuzzy and interval-valued data.

This work not only enriches the theory of fuzzy calculus on time scales but also lays a foundation for further exploration of fuzzy dynamic systems in diverse applications. Future research could focus on extending these results to broader classes of fuzzy differential and integral equations, examining their numerical solutions, and exploring real-world applications in fields such as control theory, economics, and biology. Moreover, the interplay between various types of differentiability and integrability on time scales, particularly under uncertain or imprecise conditions, remains a promising direction for further investigation.

\textbf{Author Contribution Statements:} The authors equally contributed to the paper.

\medskip

\textbf{Declaration of Competing Interests:}
The authors declare that they have no known competing financial interests or personal relationships that could have appeared to influence the work reported in this paper.


\medskip
\begin{thebibliography}{1}
\bibitem{Hilger1} Hilger, S., Ein Makettenkalkuls mit Anwendung auf Zentrumsmannigfaltigkeiten. Ph.D. Thesis, Universitat Wurzburg, Würzburg, Germany, 1988.
\bibitem{Hilger2} Hilger, S., Analysis on measure chains-A unified approach to continuous and discrete calculus, \textit{Results Math.}, 18 (1990), 18–56. \url{https://doi.org/10.1007/BF03323153}
\bibitem{Bohner} Bohner, M., Peterson, A., Dynamic Equations on Time Scales: An Introduction with Applications, Birkhauser, 2001. \url{https://doi.org/10.1007/978-1-4612-0201-1}
\bibitem{Bohner2} Bohner, M., Peterson, A. Advances in Dynamic Equations on Time Scales, Birkhauser: Boston, 2003. \url{https://doi.org/10.1007/978-0-8176-8230-9}
\bibitem{Agarwal} Agarwal, R. P., Bohner, M., Basic calculus on time scales and some of its applications, \textit{Results Math.}, 35 (1999), 3–22. \url{https://doi.org/10.1007/BF03322019}
\bibitem{Guseinov} Guseinov, G. S., Integration on time scales, \textit{Journal of Mathematical Analysis and Applications}, 285(1) (2003), 107–127. \url{https://doi.org/10.1016/S0022-247X(03)00361-5}
\bibitem{Guseinov2} Guseinov, G. S., Kaymakçalan, B.,  Basics of Riemann delta and nabla integration on time scales, \textit{Journal of Difference Equations and Applications}, 8(11) (2002), 1001-1017. \url{https://doi.org/10.1080/10236190290015272}
\bibitem{Anderson} Anderson, D., Bullock, J., Erbe, L., Peterson, A.,  Tran, H.,  Nabla dynamic equations, In Advances in dynamic equations on time scales, Boston, MA: Birkhäuser Boston, 2003. \url{https://doi.org/10.1007/978-0-8176-8230-9\_3}
\bibitem{Malinowska} Malinowska, A. B., Torres, D.F.M., On the diamond-alpha Riemann integral and mean value theorems on time scales, \textit{Dynam Systems Appl.}, 18(3–4) (2009), 469-481. \url{https://doi.org/10.48550/arXiv.0804.4420}
\bibitem{Malinowska2} Mozyrska, D., Torres, D. F. M., A study of diamond-alpha dynamic equations on regular time scales, \textit{Afr. Diaspora J. Math (NS)}, 8(1) (2009), 35-47. \url{https://doi.org/10.48550/arXiv.0902.1380}
\bibitem{Lakshmikantham2} Lakshmikantham, V., Sivasundaram, S.,  Kaymakçalan, B., Dynamic systems on measure chains, Springer Science and Business Media, 370, 2013. \url{https://doi.org/10.1007/978-1-4757-2449-3}
\bibitem{Kayar} Kayar, Z., Kaymakçalan, B., The complementary nabla Bennett-Leindler type inequalities, \textit{Communications Faculty of Sciences University of Ankara Series A1 Mathematics and Statistics}, 71(2) (2022), 349-376. \url{https://doi.org/10.31801/cfsuasmas.930138}
\bibitem{Zadeh} Zadeh, L.A., Fuzzy Sets, \textit{Information and Control}, 8 (1965), 338-353.
\bibitem{Kaleva} Kaleva, O., Fuzzy differential equations, \textit{Fuzzy Sets and Systems}, 24 (1987), 301–317. \url{https://doi.org/10.1016/0165-0114(87)90029-7}
\bibitem{Lakshmikantham} Lakshmikantham, V., Mohapatra, R. N., Theory of Fuzzy Differential Equations and Inclusions, Taylor and Francis, Abingdon, UK, 2003. \url{https://doi.org/10.1201/9780203011386}
\bibitem{Bede2} Bede, B., Gal, S. G., Generalizations of the differentiability of fuzzy-number-valued functions with application to fuzzy differential equations, \textit{Fuzzy sets and systems}, 151 (2005), 581–599. \url{https://doi.org/10.1016/j.fss.2004.08.001}
\bibitem{Bede3} Bede, B., Rudas, I. J., Bencsik, A. L., First order linear fuzzy differential equations under generalized differentiability, \textit{Information Sciences}, 177 (2007), 1648–1662. \url{https://doi.org/10.1016/j.ins.2006.08.021}
\bibitem{Li} Li, J., Zhao, A., Yan, J., Cauchy problem of fuzzy differential equations under generalized differentiability, \textit{Fuzzy Sets and Systems}, 200 (2012), 1–24. \url{https://doi.org/10.1016/j.fss.2011.10.009}
\bibitem{Stefanini2} Stefanini, L., Bede, B., Generalized Hukuhara differentiability of interval-valued functions and interval differential equations, \textit{Nonlinear Analysis: Theory, Methods and Applications}, 71 (2009), 1311–1328. \url{https://doi.org/10.1016/j.na.2008.12.005}
\bibitem{Vasavi} Vasavi, C. H., Kumar, G. S., Murty, M. S. N., Fuzzy dynamic equations on time scales under generalized delta derivative via contractive-like mapping principles, \textit{Indian Journal of Science and Technology}, 9(25) (2016), 1-6. \url{https://doi.org/10.17485/ijst/2016/v9i25/85323}
\bibitem{Vasavi1} Vasavi, C. H., Kumar, G. S., Murty, M. S. N., Fuzzy Hukuhara delta differential and applications to fuzzy dynamic equations on time scales, \textit{Journal of Uncertain Systems}, 10(3) (2016), 163-180.
\bibitem{Vasavi2} Vasavi, C. H., Kumar, G. S., Murty, M. S. N., Fuzzy dynamic equations on time scales under second type Hukuhara delta derivative, \textit{International Journal of Chemical Sciences}, 14(1) (2016), 49-66.
\bibitem{Vasavi3} Vasavi, C. H., Kumar, G. S., Murty, M. S. N., Generalized differentiability and integrability for fuzzy set-valued functions on time scales, \textit{Soft Computing}, 20(3) (2016), 1093-1104. \url{https://doi.org/10.1007/s00500-014-1569-1}
\bibitem{Fard} Fard, O. S.,  Bidgoli, T. A., Calculus of fuzzy functions on time scales (I), \textit{Soft Computing}, 19 (2015), 293-305. \url{https://doi.org/10.1007/s00500-014-1252-6}
\bibitem{you} X. X. You, D. F. Zhao, B. W. Li, Nabla-Hukuhara derivative of fuzzy valued functions on time scales, \textbf{38}(4) (2018), 580--588.
\bibitem{mert1} Mert, F. R., Bayeğ, S., Kaymakçalan, B., Generalized  Hukuhara Delta Derivatives: New Characterizations and Extensions for Fuzzy Functions on Time Scales (Submitted)
\bibitem{Sheng} Sheng, Q., Fadag, M., Henderson, J., Davis, J. M., An exploration of combined dynamic derivatives on time scales and their applications, \textit{Nonlinear Analysis: Real World Applications}, 7(3) (2006), 395-413. \url{https://doi.org/10.1016/j.nonrwa.2005.03.008}
\bibitem{Rogers} Rogers Jr, J. W., Sheng, Q., Notes on the diamond-$\alpha$ dynamic derivative on time scales, \textit{Journal of Mathematical Analysis and Applications}, 326(1)  (2007), 228-241. \url{https://doi.org/10.1016/j.jmaa.2006.03.004}
\bibitem{Truong} Truong, T., Schneider, B., Nguyen Le Toan Nhat, L., Diamond alpha differentiability of interval-valued functions and its applicability to interval differential equations on time scales, \textit{Iranian Journal of Fuzzy Systems}, 21(1) (2024), 143-158. \url{https://doi.org/10.22111/IJFS.2024.45184.7977}
\bibitem{bayeg} Bayeğ, S., Mert, F. R., Generalized Hukuhara Diamond Alpha Derivative of Fuzzy Valued Functions on Time Scales, \textit{Communications Faculty of Sciences University of Ankara Series A1 Mathematics and Statistics}, (2024), Accepted.
\bibitem{You2017} You, X., Zhao, D., Cheng, J., \& Li, T. (2017). The fuzzy C-delta integral on time scales. \textit{The Journal of Nonlinear Sciences and Applications}, 11(1), 161-171. \url{https://doi.org/10.22436/jnsa.011.01.12}
\bibitem{Leelavathi2018} Leelavathi, R., Kumar, G., \& Murty, M. (2018). Nabla integral for fuzzy functions on time scales. \textit{International Journal of Applied Mathematics}, 31(5). \url{https://doi.org/10.12732/ijam.v31i5.11}
\bibitem{Leelavathi} Leelavathi, R., Suresh Kumar, G., Agarwal, R. P., Wang, C., and Murty, M. S. N., Generalized nabla differentiability and integrability for fuzzy functions on time scales, \textit{Axioms}, 9(2) (2020), 65. \url{https://doi.org/10.3390/axioms9020065}
\bibitem{Afariogun2021} Afariogun, D., Mogbademu, A., \& Olaoluwa, H. (2021). On fuzzy Henstock-Kurzweil-Stieltjes-$\diamond$-double integral on time scales. \textit{Journal of Mathematical Analysis and Modeling}, 2(2), 38-49. \url{https://doi.org/10.48185/jmam.v2i2.295}
\bibitem{Li2024} Li, J. (2024). On fuzzy Henstock-Stieltjes integral on time scales with respect to bounded variation function. \textit{PLoS One}, 19(9), e0309031. \url{https://doi.org/10.1371/journal.pone.0309031}
\bibitem{Shahidi2020} Shahidi, M., \& Khastan, A. (2020). Linear fuzzy Volterra integral equations on time scales. \textit{Computational and Applied Mathematics}, 39(3). \url{https://doi.org/10.1007/s40314-020-01205-8}
\bibitem{Shahidi2023} Shahidi, M. (2023). A study on fuzzy Volterra integral equations for $s$-correlated fuzzy processes on time scales. \textit{Fuzzy Sets and Systems}, 471, 108695. \url{https://doi.org/10.1016/j.fss.2023.108695}

\bibitem{mert2} Mert, F. R., Bayeğ, S., Kaymakçalan, B., On the Generalized Hukuhara Nabla Differentiability of Fuzzy Functions on Time Scales via Characterization Theorem (Submitted)


\bibitem{Bourdin} Bourdin, L., Nonshifted calculus of variations on time scales with $\nabla$-differentiable $\sigma$. \textit{Journal of Mathematical Analysis and Applications}, 411(2) (2014), 543-554.


\bibitem{Negoita} Negoita, C., Ralescu, D., Application of Fuzzy Sets to System Analysis, Wiley, New York, 1975. \url{https://doi.org/10.1007/978-3-0348-5921-9}


\bibitem{Stefanini} Stefanini, L., A generalization of Hukuhara difference and division for interval and fuzzy arithmetic, \textit{Fuzzy Sets and Systems}, 161 (2010), 1564-1584. \url{https://doi.org/10.1016/j.fss.2009.06.009}


\bibitem{Bede} Bede, B., and Stefanini, L., Generalized differentiability of fuzzy-valued functions, \textit{Fuzzy sets and systems}, 230 (2013), 119-141. \url{https://doi.org/10.1016/j.fss.2012.10.003}


\bibitem{Diamond2} Diamond, P., Kloeden, P., Metric spaces of fuzzy sets, \textit{Fuzzy sets and systems}, 35(2) (1990), 241-249. \url{https://doi.org/10.1016/0165-0114(90)90197-E}

\bibitem{Dubois}
Dubois, D., Prade, H., Fuzzy numbers: an overview. In \textit{Readings in Fuzzy Sets for Intelligent Systems}, Elsevier, 1993. \url{https://doi.org/10.1016/B978-1-4832-1450-4.50015-8}

\bibitem{Goetschel}
Goetschel, R., Voxman, W., Elementary fuzzy calculus. \textit{Fuzzy Sets and Systems}, 18 (1986), 31–43. \url{https://doi.org/10.1016/0165-0114(86)90026-6}




\end{thebibliography}
\end{document}